\documentclass{amsart}

\usepackage{amssymb}
\usepackage{eepic}
\usepackage{color}

\allowdisplaybreaks

\newtheorem{theorem}{Theorem}[section]
\newtheorem{lemma}[theorem]{Lemma}

\newtheorem{corollary}[theorem]{Corollary}
\newtheorem{remark}[theorem]{Remark}

\newtheorem{Atheorem}{Theorem}

\newtheorem{Alemma}[Atheorem]{Lemma}

\newtheorem{Btheorem}{Theorem}

\newtheorem{Bremark}[Btheorem]{Remark}

\newtheorem{TheoA}{Theorem A1}
\newtheorem{TheoB}{Theorem A2}
\newtheorem{TheoC}{Theorem B1}
\newtheorem{TheoD}{Theorem B2}

\newtheorem{localtheo}{Hibert space valued pseudo-localization}

\newtheorem{Cuculescutheo}{Cuculescu's construction \cite{Cu}}
\newtheorem{GundyDecomposition}{Gundy's decomposition \cite{PR}}
\newtheorem{CZDecomposition}{Calder\'on-Zygmund decomposition \cite{P2}}
\newtheorem{shifttheo}{Shifted quasi-orthogonal decomposition}
\newtheorem{cotlar}{Cotlar lemma}
\newtheorem{schur}{Schur lemma}
\newtheorem{localest}{A localization estimate}

\newcommand{\Z}{\mathbb{Z}}
\newcommand{\R}{\mathbb{R}}
\newcommand{\C}{\mathbb{C}}

\newcommand{\summ}{\sum\nolimits}

\newcommand{\Mn}{\mathcal{A}}
\def\M{\mathcal{M}}
\def\Q{\mathcal{Q}}
\def\1{\mathbf{1}}

\newcommand{\dem}{\noindent {\bf Proof. }}

\newcommand{\fin}{\hspace*{\fill} $\square$ \vskip0.2cm}

\def\esssup{\mathop{\mathrm{ess \ sup}}}

\addtocounter{tocdepth}{-1}

\begin{document}

\title[Noncommutative Littlewood-Paley inequalities]
{Pseudo-localization of singular integrals and noncommutative
Littlewood-Paley inequalities}

\author[Tao Mei and Javier Parcet]
{Tao Mei and Javier Parcet}

\maketitle

\vskip-0.5cm

\null

\tableofcontents

\addtolength{\parskip}{+1ex}

\vskip-2cm


\section*{Introduction}

Understood in a wide sense, square functions play a central role
in classical Littlewood-Paley theory. This entails for instance
dyadic type decompositions of Fourier series, Stein's theory for
symmetric diffusion semigroups or Burkholder's martingale square
function. All these topics provide a deep technique when dealing
with quasi-orthogonality methods, sums of independent variables,
Fourier multiplier estimates... The historical survey \cite{St1}
is an excellent exposition. In a completely different setting, the
rapid development of operator space theory and quantum probability
has given rise to noncommutative analogs of several classical
results in harmonic analysis. We find new results on Fourier/Schur
multipliers, a settled theory of noncommutative martingale
inequalities, an extension for semigroups on noncommutative $L_p$
spaces of the Littlewood-Paley-Stein theory, a noncommutative
ergodic theory and a germ for a noncommutative Calder\'on-Zygmund
theory. We refer to \cite{Ha,J1,JLX,JM,JX2,P2,PX1} and the
references therein.

The aim of this paper is to produce weak type inequalities for a
large class of noncommutative square functions. In conjunction
with BMO type estimates interpolation and duality, we will obtain
the corresponding norm equivalences in the whole $L_p$ scale.
Apart from the results themselves, perhaps the main novelty relies
on our approach. Indeed, emulating the classical theory, we shall
develop a \emph{row/column} valued theory of noncommutative
martingale transforms and operator valued Calder\'on-Zygmund
operators. This seems to be new in the noncommutative setting and
may be regarded as a first step towards a noncommutative
vector-valued theory. To illustrate it, let us state our result
for noncommutative martingales.

%In what follows, $(\mathcal{M}_n)_{n \ge 1}$ will stand for a
%weak$^*$ dense increasing filtration in a semifinite von Neumann
%algebra $(\mathcal{M},\tau)$ equipped with a normal semifinite
%faithful trace $\tau$. Our weak type inequality reads as follows.

\begin{TheoA} Let $(\mathcal{M}_n)_{n \ge 1}$ stand for a
weak$^*$ dense increasing filtration in a semifinite von Neumann
algebra $(\mathcal{M},\tau)$ equipped with a normal semifinite
faithful trace $\tau$. Given $f = (f_n)_{n \ge 1}$ an $L_1(\M)$
martingale, let
$$T_m \hskip-1pt f = \sum_{k=1}^\infty \xi_{km} df_k \quad \mbox{with}
\quad \sup_{k \ge 1} \sum_{m=1}^\infty |\xi_{km}|^2 \lesssim 1.$$
Then, there exists a decomposition $T_m \hskip-1pt f = A_m
\hskip-1pt f + B_m \hskip-1pt f$, satisfying
$$\Big\| \Big( \sum_{m=1}^\infty (A_m \hskip-1pt f)(A_m \hskip-1pt
f)^* \Big)^\frac12 \Big\|_{1,\infty} + \Big\| \Big(
\sum_{m=1}^\infty (B_m \hskip-1pt f)^*(B_m \hskip-1pt f)
\Big)^\frac12 \Big\|_{1,\infty} \lesssim \, \sup_{n \ge 1}
\|f_n\|_1.$$
\end{TheoA}

In the statement above, $df_k$ denotes the $k$-th martingale
difference of $f$ relative to the filtration $(\mathcal{M}_n)_{n
\ge 1}$ and $\| \ \|_{1,\infty}$ refers to the norm on
$L_{1,\infty}(\mathcal{M})$. In the result below, we also need to
use the norm on $\mathrm{BMO}(\mathcal{M})$ relative to our
filtration as well as the norm on $L_p(\mathcal{M}; \ell^2_{r \!
c})$. All these norms are standard in the noncommutative setting
and we refer to Section \ref{S1} below for precise definitions.
Moreover, in what follows $\delta_k$ and $e_{ij}$ will stand for
unit vectors of sequence spaces and matrix algebras respectively.

\begin{TheoB}
Let us set $\mathcal{R} = \M \bar\otimes \mathcal{B}(\ell_2)$.
Assume that $f$ is an $L_\infty(\mathcal{M})$ martingale relative
to the filtration $(\M_n)_{n \ge 1}$ and define $T_m \hskip-1pt f$
with coefficients $\xi_{km}$ satisfying the same condition above.
Then, we have
$$\Big\| \sum_{m=1}^\infty T_m \hskip-1pt f \otimes e_{1m}
\Big\|_{\mathrm{BMO}({\mathcal{R}})} + \Big\| \sum_{m=1}^\infty
T_m \hskip-1pt f \otimes e_{m1}
\Big\|_{\mathrm{BMO}({\mathcal{R}})} \lesssim \, \sup_{n \ge 1}
\|f_n\|_\infty.$$ Therefore, given $1 < p < \infty$ and $f \in
L_p(\mathcal{M})$, we deduce $$\Big\| \sum_{m=1}^\infty T_m
\hskip-1pt f \otimes \delta_m \Big\|_{L_p(\mathcal{M}; \ell^2_{r
\! c})} \le \, c_p \, \|f\|_p.$$ Moreover, the reverse inequality
also holds if $\sum_m |\xi_{km}|^2 \sim 1$ uniformly on $k$.
\end{TheoB}

Let us briefly analyze Theorems A1 and A2. Taking $\xi_{km}$ to be
the Dirac delta on $(k,m)$, we find $T_m \hskip-1pt f=df_m$ and
our results follow from the noncommutative Burkholder-Gundy
inequalities \cite{PX1,R3}. Moreover, taking $\xi_{km}=0$ for
$m>1$ we simply obtain a martingale transform with scalar
coefficients and our results follow from \cite{R1}. Other known
examples appear by considering $(\xi_{km})$ of \emph{diagonal-like
shape}. For instance, taking an arbitrary partition
$$\mathbb{N} = \bigcup_{m \ge 1} \Omega_m \quad \mbox{and} \quad
\xi_{km} = \begin{cases} 1 & \mbox{if} \ k \in \Omega_m, \\ 0 &
\mbox{otherwise}. \end{cases}$$ It is apparent that $\sum_m
|\xi_{km}|^2 = 1$ and Theorem A2 gives e.g. for $2 \le p < \infty$
$$\Big\| \Big( \sum_{m=1}^\infty \Big| \sum_{k \in \Omega_m} df_k
\, \Big|^2 \, \Big)^\frac12 \Big\|_p + \Big\| \Big(
\sum_{m=1}^\infty \Big| \sum_{k \in \Omega_m} df_k^* \, \Big|^2 \,
\Big)^\frac12 \Big\|_p \sim c_p \, \|f\|_p.$$ Except for $p=1$,
this follows from the noncommutative Khintchine inequality in
conjunction with the $L_p$ boundedness of martingale transforms.
The new examples appear when considering more general matrices
$(\xi_{km})$ and will be further analyzed in the body of the
paper.

In the framework of Theorems A1 and A2, the arguments in
\cite{PX1,R1,R3} are no longer valid. Instead, we think in our
square functions as martingale transforms with row/column valued
coefficients
\begin{eqnarray*}
\Big( \sum_{m=1}^\infty (T_m \hskip-1pt f)(T_m \hskip-1pt f)^*
\Big)^\frac12 & \sim & \sum_{k=1}^\infty \big( df_k \otimes e_{1,1} \big)\,\Big(
\overbrace{\sum_{m=1}^\infty \xi_{k,m} \mathbf{1}_\mathcal{M}
\otimes e_{1m}}^{\xi_k^r} \Big), \\
\Big( \sum_{m=1}^\infty (T_m \hskip-1pt f)^*(T_m \hskip-1pt f)
\Big)^\frac12 & \sim & \sum_{k=1}^\infty \Big(
\underbrace{\sum_{m=1}^\infty \xi_{k,m} \mathbf{1}_\mathcal{M}
\otimes e_{m1}}_{\xi_k^c} \Big) \, \big( df_k \otimes e_{1,1}
\big),
\end{eqnarray*}
where $\sim$ means to have the same $L_{1,\infty}(\mathcal{M})$ or
$L_p(\mathcal{M})$ norm. Tensorizing with the identity on
$\mathcal{B}(\ell_2)$, we have $df_k \otimes e_{1,1}=d(f \otimes
e_{1,1})_k$ and we find our row/column valued transforms.
According to \cite{R1}, we might expect
$$\Big\| \sum_{k=1}^\infty \xi_k^r \, d(f \otimes e_{1,1})_k
\Big\|_p \le \, c_p \, \sup_{k \ge 1}
\|\xi_k^r\|_{\mathcal{B}(\ell_2)} \, \Big\| \sum_{k=1}^\infty df_k
\otimes e_{1,1} \Big\|_p \lesssim \, c_p \, \|f\|_p$$ and the same
estimate for the $\xi_k^c$'s. However, it is essential in
\cite{R1} to have commuting coefficients $\xi_k \in
\mathcal{R}_{k-1} \cap \mathcal{R}_k'$, where $\mathcal{R}_n =
\mathcal{M}_n \bar\otimes \mathcal{B}(\ell_2)$ in our setting.
This is not the case. In fact, the inequality above is false in
general (e.g. take again $\xi_{km}=\delta_{(k,m)}$ with $1 < p <
2$) and Theorems A1 and A2 might be regarded as the right
substitute. The same phenomenon will occur is the context of
operator-valued Calder\'on-Zygmund operators below.

Our main tools to overcome it will be the noncommutative forms of
Gundy's and Calder\'on-Zygmund decompositions \cite{P2,PR} for
martingales transforms and singular integral operators
respectively. As it was justified in \cite{P2}, there exists
nevertheless a substantial difference between both settings.
Namely, martingale transforms are local operators while
Calder\'on-Zygmund operators are only pseudo-local. In this paper
we will illustrate this point by means of Rota's dilation theorem
\cite{Ro}. The pseudo-localization estimate that we need in this
setting, to pass from martingale transforms to Calder\'on-Zygmund
operators, is a Hilbert space valued version of that given in
\cite{P2} and will be sketched in Appendix A.

Now we formulate our results for Calder\'on-Zygmund operators. Let
$\Delta$ denote the diagonal of $\R^n \times \R^n$ and fix a
Hilbert space $\mathcal{H}$. We will write in what follows $T$ to
denote an integral operator associated to a kernel $k: \R^{2n}
\setminus \Delta \to \mathcal{H}$. This means that for any smooth
test function $f$ with compact support, we have
$$Tf(x) = \int_{\R^n} k(x,y) f(y) \, dy \quad \mbox{for all} \quad
x \notin \mbox{supp} \hskip1pt f.$$ Given two points $x,y \in
\R^n$, the distance $|x-y|$ between $x$ and $y$ will be taken for
convenience with respect to the $\ell_\infty(n)$ metric. As usual,
we impose size and smoothness conditions on the kernel:
\begin{itemize}
\item[a)] If $x,y \in \R^n$, we have $$\big\| k(x,y)
\big\|_\mathcal{H} \ \lesssim \ \frac{1}{|x-y|^n}.$$
\item[b)] There exists $0 < \gamma \le 1$ such that
$$\begin{array}{rcl} \big\| k(x,y) - k(x',y) \big\|_\mathcal{H}
& \lesssim & \displaystyle \frac{|x-x'|^\gamma}{|x-y|^{n+\gamma}}
\quad \mbox{if} \quad |x-x'| \le \frac12 \hskip1pt |x-y|, \\
[10pt] \big\| k(x,y) - k(x,y') \big\|_\mathcal{H} & \lesssim &
\displaystyle \frac{|y-y'|^\gamma}{|x-y|^{n+\gamma}} \quad
\mbox{if} \quad |y-y'| \le \frac12 \hskip1pt |x-y|.
\end{array}$$
\end{itemize}
We will refer to this $\gamma$ as the Lipschitz parameter of the
kernel. The statement of our results below requires to consider
appropriate $\mathcal{H}$-valued noncommutative function spaces as
in \cite{JLX}. Let us first consider the algebra $\Mn_B$ of
essentially bounded functions with values in $\M$
$$\Mn_B = \Big\{ f: \R^n \to \M \, \big| \, f \ \mbox{strongly
measurable s.t.} \ \esssup_{x \in \R^n} \|f(x)\|_{\M} < \infty
\Big\},$$ equipped with the n.s.f. trace $\varphi(f) = \int_{\R^n}
\tau(f(x)) \, dx$. The weak-operator closure $\Mn$ of $\Mn_B$ is a
von Neumann algebra. Given a norm 1 element $e \in \mathcal{H}$,
take $p_e$ to be the orthogonal projection onto the
one-dimensional subspace generated by $e$ and define
\begin{eqnarray*}
L_p(\Mn; \mathcal{H}_r) & = & (\mathbf{1}_\Mn \otimes
p_e) L_p(\Mn \bar\otimes \mathcal{B}(\mathcal{H})), \\
L_p(\Mn; \mathcal{H}_c) & = & L_p(\Mn \bar\otimes
\mathcal{B}(\mathcal{H})) (\mathbf{1}_\Mn \otimes p_e).
\end{eqnarray*}
This definition is essentially independent of the choice of $e$.
Indeed, given a function $f \in L_p(\Mn; \mathcal{H}_r)$ we may
regard it as an element of $L_p(\Mn \bar\otimes
\mathcal{B}(\mathcal{H}))$, so that the product $ff^*$ belongs to
$(\mathbf{1}_\Mn \otimes p_e) L_{p/2}(\Mn \bar\otimes
\mathcal{B}(\mathcal{H})) (\mathbf{1}_\Mn \otimes p_e)$ which may
be identified with $L_{p/2}(\Mn)$. When $f \in L_p(\Mn;
\mathcal{H}_c)$ the same holds for $f^*f$ and we conclude
$$\|f\|_{L_p(\Mn; \mathcal{H}_r)} = \big\| (ff^*)^\frac12
\big\|_{L_p(\Mn)} \quad \mbox{and} \quad \|f\|_{L_p(\Mn;
\mathcal{H}_c)} = \big\| (f^*f)^\frac12 \big\|_{L_p(\Mn)}.$$
Arguing as in \cite[Chapter 2]{JLX}, we may use these identities
to regard $L_p(\Mn) \otimes \mathcal{H}$ as a dense subspace of
$L_p(\Mn;\mathcal{H}_r)$ and $L_p(\Mn;\mathcal{H}_c)$. More
specifically, given a function $f = \sum_k g_k \otimes v_k \in
L_p(\Mn) \otimes \mathcal{H}$, we have
\begin{eqnarray*}
\|f\|_{L_p(\Mn;\mathcal{H}_r)} & = & \Big\| \Big( \sum_{i,j}
\langle v_i,v_j \rangle \, g_ig_j^* \Big)^\frac12
\Big\|_{L_p(\Mn)}, \\ \|f\|_{L_p(\Mn;\mathcal{H}_c)} & = & \Big\|
\Big( \sum_{i,j} \langle v_i,v_j \rangle \, g_i^*g_j \Big)^\frac12
\Big\|_{L_p(\Mn)}.
\end{eqnarray*}
This procedure may also be used to define the spaces
$$L_{1,\infty}(\mathcal{A}; \mathcal{H}_r) \quad \mbox{and} \quad
L_{1,\infty} (\Mn; \mathcal{H}_c).$$ It is clear that
$L_2(\mathcal{M}; \mathcal{H}_{r}) = L_2(\mathcal{M};
\mathcal{H}_{c})$ and we will denote it by $L_2(\mathcal{M};
\mathcal{H}_{oh})$.

\begin{TheoC} Given $f \in L_1(\mathcal{A})$, define formally
$$T \! f(x) = \int_{\R^n} k(x,y) f(y) \, dy$$ where the kernel
$k: \R^{2n} \setminus \Delta \to \mathcal{H}$ satisfies the
size/smoothness conditions imposed above. Assume further that $T$
defines a bounded map $L_2(\mathcal{A}) \to L_2(\mathcal{A};
\mathcal{H}_{oh})$. Then we may find a decomposition $T \! f = A f
+ B \hskip-1pt f$, satisfying $$\big\| Af
\big\|_{L_{1,\infty}(\Mn; \mathcal{H}_r)} + \big\| B \hskip-1pt f
\big\|_{L_{1,\infty}(\Mn; \mathcal{H}_c)} \lesssim \, \|f\|_1.$$
\end{TheoC}

To state the following result, we need to define the corresponding
$\mathcal{H}$-valued $\mathrm{BMO}$ norm. Assume for simplicity
that $\mathcal{H}$ is separable and fix an orthonormal basis
$(v_k)_{k \ge 1}$ in $\mathcal{H}$. Then, given a norm 1 element
$e \in \mathcal{H}$, we may identify (as above) $f = \sum_k g_k
\otimes v_k \in \mathrm{BMO}(\Mn) \otimes \mathcal{H}$ with
$$\begin{array}{rcccl} {}_ef & = & \summ_k g_k \otimes (e \otimes
v_k) & = & (\mathbf{1}_\Mn \otimes p_e) \Big( \summ_k g_k \otimes
(e \otimes v_k) \Big), \\ f_e & = & \summ_k g_k \otimes (v_k
\otimes e) & = & \Big( \summ_k g_k \otimes (e \otimes v_k) \Big)
(\mathbf{1}_\Mn \otimes p_e),
\end{array}$$ where $e \otimes v_k$ is understood as the rank 1
operator $\xi \in \mathcal{H} \mapsto \langle v_k, \xi \rangle e$
and $v_k \otimes e$ stands for $\xi \in \mathcal{H} \mapsto
\langle e,\xi \rangle v_k$. Then we define the spaces
$\mathrm{BMO}(\Mn;\mathcal{H}_r)$ and
$\mathrm{BMO}(\Mn;\mathcal{H}_c)$ as the closure of
$\mathrm{BMO}(\Mn) \otimes \mathcal{H}$ with respect to the norms
$$\|f\|_{\mathrm{BMO}(\Mn;\mathcal{H}_r)} =
\|{}_ef\|_{\mathrm{BMO}(\Mn \otimes \mathcal{B}(\mathcal{H}))}
\quad \mbox{and} \quad \|f\|_{\mathrm{BMO}(\Mn;\mathcal{H}_c)} =
\|f_e\|_{\mathrm{BMO}(\Mn \otimes \mathcal{B}(\mathcal{H}))}.$$ In
the following result, we also use the standard terminology
$$L_p(\mathcal{A}; \mathcal{H}_{r \! c}) = \left\{ \begin{array}{ll}
L_p(\mathcal{A}; \mathcal{H}_r) + L_p(\mathcal{A}; \mathcal{H}_c)
& 1 \le p \le 2, \\ L_p(\mathcal{A}; \mathcal{H}_r) \cap \hskip1pt
L_p(\mathcal{A}; \mathcal{H}_c) & 2 \le p \le \infty.
\end{array} \right.$$

\begin{TheoD}
If $f \in L_\infty(\mathcal{A})$, we also have
$$\big\| T \! f \big\|_{\mathrm{BMO}(\Mn;\mathcal{H}_r)} +
\big\| T \! f \big\|_{\mathrm{BMO}(\Mn;\mathcal{H}_c)} \lesssim \,
\|f\|_\infty.$$ Therefore, given $1 < p < \infty$ and $f \in
L_p(\mathcal{A})$, we deduce $$\big\| T \! f
\big\|_{L_p(\mathcal{A}; \mathcal{H}_{r \! c})} \le \, c_p \,
\|f\|_p.$$ Moreover, the reverse inequality holds whenever $\|T \!
f\|_{L_2(\Mn; \mathcal{H}_{oh})} = \|f\|_{L_2(\Mn)}$.
\end{TheoD}

%Note that $L_p(\mathcal{A}; \mathcal{H}_{r \! c}) = L_p(\M;
%\ell_{r \! c}^2)$ for $(\Mn, \mathcal{H}) = (\M, \ell_2)$ and that
%we might formulate Theorems A1 and A2 for $\mathcal{H}$-valued
%coefficients, but we do not have in mind interesting examples.

In Section \ref{S1} we prove Theorems A1 and A2. Then we study an
specific example on ergodic averages as in \cite{St2}. In
conjunction with Rota's theorem, this shows the relevance of
pseudo-localization in the Calder\'on-Zygmund setting. We also
find some multilinear and operator-valued forms of our results.
Theorems B1 and B2 are proved in Section \ref{S2}. The proof
requires a Hilbert space valued pseudo-localization estimate
adapted from \cite{P2} in Appendix A. After the proof, we list
some examples and applications. Although most of the examples are
\emph{semicommutative}, we find several new square functions not
considered in \cite{JLX,Me} and find an application in the
\emph{fully} noncommutative setting which will be explored in
\cite{JMP}. Finally, following a referee's suggestion, we also
include an additional Appendix B with some background on
noncommutative $L_p$ spaces, noncommutative martingales and a few
examples for nonexpert readers.

\numberwithin{equation}{section}

\section{Martingale transforms}
\label{S1}

In this section, we prove Theorems A1 and A2. As a preliminary, we
recall the definition of some noncommutative function spaces and
the statement of some auxiliary results. We shall assume that the
reader is familiar with noncommutative $L_p$ spaces. Given
$(\M,\tau)$ a semifinite von Neumann algebra equipped with a
n.s.f. trace, the noncommutative weak $L_1$-space
$L_{1,\infty}(\mathcal{M})$ is defined as the set of all
$\tau$-measurable operators $f$ for which the quasi-norm
$$\left\|f\right\|_{1,\infty} = \sup_{\lambda > 0} \, \lambda
\hskip1pt \tau \Big\{ |f| > \lambda \Big\}$$ is finite. In this
case, we write $\tau \big\{ |f| > \lambda \big\}$ to denote the
trace of the spectral projection of $|f|$ associated to the
interval $(\lambda,\infty)$. We find this terminology more
intuitive, since it is reminiscent of the classical one. The space
$L_{1,\infty}(\M)$ satisfies a quasi-triangle inequality that will
be used below with no further reference $$\lambda \, \tau \Big\{
|f_1+f_2| > \lambda \Big\} \le \lambda \, \tau \Big\{ |f_1| >
\lambda/2 \Big\} + \lambda \, \tau \Big\{ |f_2| > \lambda/2
\Big\}.$$ We refer the reader to \cite{FK,PX2} for a more in depth
discussion on these notions.

Let us now define the space $\mathrm{BMO}(\M)$. Let $L_0(\M)$
stand for the $*$-algebra of $\tau$-measurable operators
affiliated to $\M$ and fix a filtration $(\M_n)_{n \ge 1}$. Let us
write $\mathsf{E}_n: \M \to \M_n$ for the corresponding
conditional expectation. Then we define $\mathrm{BMO}_\M^r$ and
$\mathrm{BMO}_\M^c$ as the spaces of operators $f \in L_0(\M)$
with norm (modulo multiples of $\mathbf{1}_\M$)
\begin{eqnarray*}
\|f\|_{\mathrm{BMO}_\M^r} & = & \sup_{n \ge 1} \Big\| \mathsf{E}_n
\Big( (f - \mathsf{E}_{n-1}(f))(f - \mathsf{E}_{n-1}(f))^*
\Big)^\frac12 \Big\|_\M, \\ \|f\|_{\mathrm{BMO}_\M^c} & = &
\sup_{n \ge 1} \Big\| \mathsf{E}_n \Big( (f -
\mathsf{E}_{n-1}(f))^*(f - \mathsf{E}_{n-1}(f)) \Big)^\frac12
\Big\|_\M.
\end{eqnarray*}
It is easily checked that we have the identities
\begin{eqnarray*}
\|f\|_{\mathrm{BMO}_\M^r} & = & \sup_{n \ge 1} \Big\| \mathsf{E}_n
\Big( \sum_{k \ge n} df_k df_k^* \Big)^\frac12 \Big\|_\M, \\
\|f\|_{\mathrm{BMO}_\M^c} & = & \sup_{n \ge 1} \Big\| \mathsf{E}_n
\Big( \sum_{k \ge n} df_k^* df_k \Big)^\frac12 \Big\|_\M.
\end{eqnarray*}
We define $\mathrm{BMO}(\M) = \mathrm{BMO}_\M^r \cap
\mathrm{BMO}_\M^c$ with norm given by $$\|f\|_{\mathrm{BMO}(\M)} =
\max \Big\{ \|f\|_{\mathrm{BMO}_\M^r}, \|f\|_{\mathrm{BMO}_\M^c}
\Big\}.$$ Finally, the space $L_p(\M; \ell_{r \! c}^2)$ was
already defined in the Introduction.

A key tool in proving weak type inequalities for noncommutative
martingales is due to Cuculescu. It can be viewed as a
noncommutative analogue of the weak type $(1,1)$ boundedness of
Doob's maximal function.

\begin{Cuculescutheo}
Let $f = (f_1, f_2, \ldots)$ be a positive $L_1$ martingale
relative to the filtration $(\mathcal{M}_n)_{n \ge 1}$ and let
$\lambda$ be a positive number. Then there exists a decreasing
sequence of projections
$$q(\lambda)_1, q(\lambda)_2, q(\lambda)_3, \ldots$$ in
$\mathcal{M}$ satisfying the following properties
\begin{itemize}
\item[i)] $q(\lambda)_n$ commutes with $q(\lambda)_{n-1} f_n
q(\lambda)_{n-1}$ for each $n \ge 1$.
\item[ii)] $q(\lambda)_n$ belongs to $\mathcal{M}_n$ for each $n
\ge 1$ and $q(\lambda)_n f_n q(\lambda)_n \le \lambda \hskip1pt
q(\lambda)_n$.
\item[iii)] The following estimate holds $$\tau \Big(
\mathbf{1}_\M - \bigwedge_{n \ge 1} q(\lambda)_n \Big) \le
\frac{1}{\lambda} \hskip1pt \sup_{n \ge 1} \|f_n\|_1.$$
\end{itemize}
Explicitly, we set $q(\lambda)_0 = \1_\M$ and define $q(\lambda)_n
= \chi_{(0,\lambda]}(q(\lambda)_{n-1} f_n q(\lambda)_{n-1})$.
\end{Cuculescutheo}

Another key tool for what follows is Gundy's decomposition for
noncommutative martingales. We need a weak notion of support which
is quite useful when dealing with weak type inequalities. For a
non-necessarily self-adjoint $f \in \M$, the two sided null
projection of $f$ is the greatest projection $q$ in $\M$
satisfying $qfq = 0$. Then we define the \emph{weak support
projection} of $f$ as $$\mathrm{supp}^* f = \1_\Mn - q.$$ It is
clear that $\mathrm{supp}^* f = \mathrm{supp} f$ when $\M$ is
abelian. Moreover, this notion is weaker than the usual support
projection in the sense that we have $\mathrm{supp}^* f \le
\mathrm{supp} f$ for any self-adjoint $f \in \M$ and
$\mathrm{supp}^* f$ is a subprojection of both the left and right
supports in the non-self-adjoint case.

\begin{GundyDecomposition} Let $f = (f_1, f_2,
\ldots)$ be a positive $L_1$ martingale relative to the filtration
$(\mathcal{M}_n)_{n \ge 1}$ and let $\lambda$ be a positive
number. Then $f$ can be decomposed $f = \alpha+\beta+\gamma$ as
the sum of three martingales relative to the same filtration and
satisfying $$\max \Big\{\frac{1}{\lambda} \, \sup_{n \ge 1}
\|\alpha_n\|_2^2, \, \sum_{k=1}^{\infty} \|d \beta_k\|_1, \,
\lambda \tau \Big( \bigvee_{k \ge 1} \mathrm{supp}^* d \gamma_k
\Big) \Big\} \lesssim \, \sup_{n \ge 1} \|f_n\|_1.$$ We may write
$\alpha, \beta$ and $\gamma$ in terms of their martingale
differences
\begin{eqnarray*}
d\alpha_k & = & q_k(\lambda) df_k q_k(\lambda) - \mathsf{E}_{k-1}
\big( q_k(\lambda) df_k q_k(\lambda) \big), \\ d\beta_k & = &
q_{k-1}(\lambda) df_k q_{k-1}(\lambda) - q_k(\lambda) df_k
q_k(\lambda) + \mathsf{E}_{k-1} \big( q_k(\lambda) df_k
q_k(\lambda) \big), \\ d\gamma_k & = & df_k - q_{k-1}(\lambda)
df_k q_{k-1}(\lambda).
\end{eqnarray*}
\end{GundyDecomposition}

\subsection{Weak type $(1,1)$ boundedness}

Here we prove Theorem A1. Let us begin with some harmless
assumptions. First, we shall assume that $\M$ is a finite von
Neumann algebra with a normalized trace $\tau$. The passage to the
semifinite case is just technical. Moreover, we shall sketch it in
Section \ref{S2} since the von Neumann algebra $\Mn$ we shall work
with can not be finite. Second, we may assume that the martingale
$f$ is positive and finite, so that we may use Cuculescu's
construction and Gundy's decomposition for $f$ and moreover we do
not have to worry about convergence issues.

Now we provide the decomposition $T_m \hskip-1pt f = A_m
\hskip-1pt f + B_m \hskip-1pt f$. If $(q_n(\lambda))_{n \ge 1}$
denotes the Cuculescu's projections associated to $(f,\lambda)$,
let us write in what follows $q(\lambda)$ for the projection
$$q(\lambda) = \bigwedge_{n \ge 1} q_n(\lambda).$$ Then we define
the projections $$\pi_0 = \bigwedge_{s \ge 0} q(2^s) \quad
\mbox{and} \quad \pi_k = \bigwedge_{s \ge k} q(2^s) - \bigwedge_{s
\ge k-1} q(2^s)$$ for $k \ge 1$. Since $\sum_{k \ge 0} \pi_k =
\mathbf{1}_\M$, we may write $$df_k = \sum_{i \ge j} \pi_i df_k
\pi_j + \sum_{i < j} \pi_i df_k \pi_j = \Delta_r(df_k) +
\Delta_c(df_k).$$ Then, our decomposition $T_m \hskip-1pt f = A_m
\hskip-1pt f + B_m \hskip-1pt f$ is given by $$A_m \hskip-1pt f =
\sum_{k=1}^\infty \xi_{km} \Delta_r(df_k) \quad \mbox{and} \quad
B_m \hskip-1pt f = \sum_{k=1}^\infty \xi_{km} \Delta_c(df_k).$$
Since both terms can be handled in a similar way, we shall only
prove that $$\sup_{\lambda > 0} \lambda \, \tau \Big\{ \Big(
\sum_{m=1}^\infty (A_m \hskip-1pt f) (A_m \hskip-1pt f)^*
\Big)^\frac12 > \lambda \Big\} \lesssim \, \sup_{n \ge 1}
\|f\|_1.$$ By homogeneity, we may assume that the right hand side
equals 1. This means in particular that we may also assume that
$\lambda \ge 1$ since, by the finiteness of $\M$, the left hand
side is bounded above by $1$ for $0 < \lambda < 1$. Moreover, up
to a constant $2$ it suffices to prove the result for $\lambda$
being a nonnegative power of $2$. Let us fix a nonnegative integer
$\ell$, so that $\lambda = 2^\ell$ for the rest of the proof.

\noindent Let us define $$w_\ell = \bigwedge_{s \ge \ell}
q(2^s).$$ By the quasi-triangle inequality, we are reduced to
estimate $$\lambda \, \tau \Big\{ w_\ell \Big( \sum_{m=1}^\infty
(A_m \hskip-1pt f) (A_m \hskip-1pt f)^* \Big) w_\ell > \lambda^2
\Big\} + \lambda \, \tau (\mathbf{1}_\M - w_\ell) = \mathsf{A}_1 +
\mathsf{A}_2.$$ According to Cuculescu's theorem, $\mathsf{A}_2$
is dominated by $$ \lambda \, \sum_{s \ge \ell} \tau(\mathbf{1}_\M
- q(2^s)) \le 2^\ell \Big( \sum_{s \ge \ell} \frac{1}{2^s} \Big)
\, \sup_{n \ge 1} \, \|f_n\|_1 \le 2 \, \sup_{n \ge 1} \,
\|f_n\|_1.$$ Let us now proceed with the term $\mathsf{A}_1$. We
first notice that $w_\ell \pi_k = \pi_k w_\ell = 0$ for any
integer $k > \ell$. Therefore, we find $w_\ell \Delta_r(df_k) =
w_\ell (\sum_{j \le i \le \ell} \pi_i df_k \pi_j) = w_\ell
\Delta_{r \ell}(df_k)$. Similarly, we have $\Delta_r(df_k)^*
w_\ell = \Delta_{r \ell}(df_k)^* w_\ell$ and letting $$A_{m\ell}
\hskip-1pt f = \sum_{k=1}^\infty \xi_{km} \Delta_{r\ell}(df_k),$$
we conclude $$\mathsf{A}_1 = \lambda \, \tau \Big\{ w_\ell \Big(
\sum_{m=1}^\infty (A_{m \ell} \hskip-1pt f) (A_{m \ell} \hskip-1pt
f)^* \Big) w_\ell > \lambda^2 \Big\}.$$ Moreover, using the fact
that the spectral projections $\chi_{(\lambda,\infty)}(xx^*)$ and
$\chi_{(\lambda,\infty)}(x^*x)$ are Murray-von Neumann equivalent,
we may kill the projection $w_\ell$ above and obtain the
inequality $\mathsf{A}_1 \le \lambda \, \tau \big\{ \sum_m (A_{m
\ell} \hskip-1pt f) (A_{m \ell} \hskip-1pt f)^* > \lambda^2
\big\}$. Now we use Gundy's decomposition for $(f,\lambda)$ and
quasi-triangle inequality to get
\begin{eqnarray*}
\mathsf{A}_1 & \lesssim & \lambda \, \tau \Big\{ \sum_{m=1}^\infty
(A_{m \ell} \alpha) (A_{m \ell} \alpha)^* > \lambda^2 \Big\}
\\ & + & \lambda \, \tau \Big\{
\sum_{m=1}^\infty (A_{m \ell} \beta) (A_{m \ell} \beta)^* >
\lambda^2 \Big\} \\ & + & \lambda \, \tau \Big\{ \sum_{m=1}^\infty
(A_{m \ell} \gamma) (A_{m \ell} \gamma)^* > \lambda^2 \Big\} \ = \
\mathsf{A}_\alpha + \mathsf{A}_\beta + \mathsf{A}_\gamma.
\end{eqnarray*}
We claim that $\mathsf{A}_\gamma$ is identically $0$. Indeed, note
that $$\Delta_{r \ell}(d\gamma_k) = \sum_{j \le i \le \ell} \pi_i
\Big( df_k - q_{k-1}(2^\ell) df_k q_{k-1}(2^\ell) \Big) \pi_j =
0$$ since $\pi_i q_{k-1}(2^\ell) = \pi_i$ and $q_{k-1}(2^\ell)
\pi_j = \pi_j$ for $i,j \le \ell$. Therefore, it remains to
control the terms $\mathsf{A}_\alpha$ and $\mathsf{A}_\beta$. Let
us begin with $\mathsf{A}_\alpha$. Applying Fubini to the sum
defining it, we obtain $$\sum_{m=1}^\infty (A_{m \ell} \alpha)
(A_{m \ell} \alpha)^* = \sum_{j,k=1}^\infty \Big(
\sum_{m=1}^\infty \xi_{jm} \overline{\xi_{km}} \Big)
\Delta_{r\ell}(d\alpha_j) \Delta_{r\ell}(d\alpha_k)^*.$$ It will
be more convenient, to write this as follows
\begin{eqnarray*}
\lefteqn{\sum_{m=1}^\infty (A_{m \ell} \alpha) (A_{m \ell}
\alpha)^*} \\ & = & \Big( \sum_{j=1}^\infty
\Delta_{r\ell}(d\alpha_j) e_{1j} \Big) \Big( \sum_{j,k=1}^\infty
\Big[ \sum_{m=1}^\infty \xi_{jm} \overline{\xi_{km}} \Big] e_{jk}
\Big) \Big( \sum_{k=1}^\infty \Delta_{r\ell}(d\alpha_k)^* e_{k1}
\Big) \\ & = & \Big( \sum_{j=1}^\infty \Delta_{r\ell}(d\alpha_j)
e_{1j} \Big) \Big( \sum_{j,k=1}^\infty \xi_{jk} e_{jk} \Big) \Big(
\sum_{j,k=1}^\infty \xi_{jk} e_{jk} \Big)^* \Big(
\sum_{k=1}^\infty \Delta_{r\ell}(d\alpha_k)^* e_{k1} \Big).
\end{eqnarray*}
In particular, Chebychev's inequality gives
\begin{eqnarray*}
\mathsf{A}_\alpha & \le & \frac{1}{\lambda} \, \Big\| \Big(
\sum_{j=1}^\infty \Delta_{r\ell}(d\alpha_j) e_{1j} \Big) \hskip2pt
\Big( \sum_{j,k=1}^\infty \xi_{jk} e_{jk} \Big) \Big\|_{L_2(\M
\bar\otimes \mathcal{B}(\ell_2))}^2 \\ & = & \frac{1}{\lambda} \,
\Big\| \widehat{\Delta_{r\ell}} \Big[ \Big( \sum_{j=1}^\infty
d\alpha_j e_{1j} \Big) \Big( \sum_{j,k=1}^\infty \xi_{jk} e_{jk}
\Big) \Big] \Big\|_{L_2(\M \bar\otimes \mathcal{B}(\ell_2))}^2,
\end{eqnarray*}
with $\widehat{\Delta_{r\ell}} = \Delta_{r\ell} \otimes
id_{\mathcal{B}(\ell_2)}$ a triangular truncation, bounded on
$L_2(\M \bar\otimes \mathcal{B}(\ell_2))$. Thus, we get
\begin{eqnarray*}
\mathsf{A}_\alpha & \lesssim & \frac{1}{\lambda} \, \Big\| \Big(
\sum_{j=1}^\infty d\alpha_j e_{1j} \Big) \Big( \sum_{j,k=1}^\infty
\xi_{jk} e_{jk} \Big) \Big\|_{L_2(\M \bar\otimes
\mathcal{B}(\ell_2))}^2 \\ & = & \frac{1}{\lambda} \,
\sum_{j,k=1}^\infty \Big( \sum_{m=1}^\infty \xi_{jm}
\overline{\xi_{km}} \Big) \tau(d\alpha_j d\alpha_k^*) \ = \ \Big(
\sup_{k \ge 1} \sum_{m=1}^\infty |\xi_{km}|^2 \Big) \,
\frac{1}{\lambda} \, \sum_{k=1}^\infty \tau(d\alpha_k
d\alpha_k^*).
\end{eqnarray*}
Therefore, the estimate for $\mathsf{A}_\alpha$ follows from our
hypothesis on the $\xi_{km}$'s and from the estimate for the
$\alpha$-term in Gundy's decomposition. Let us finally estimate
the term $\mathsf{A}_\beta$. Arguing as above, we clearly have
\begin{eqnarray*}
\mathsf{A}_\beta & \le & \Big\| \Big( \sum_{m=1}^\infty (A_{m
\ell} \beta) (A_{m \ell} \beta)^* \Big)^\frac12 \Big\|_{1,\infty}
\\ & = & \Big\| \widehat{\Delta_{r\ell}} \Big[ \Big( \sum_{j=1}^\infty
d\beta_j e_{1j} \Big) \Big( \sum_{j,k=1}^\infty \xi_{jk} e_{jk}
\Big) \Big] \Big\|_{L_{1,\infty}(\M \bar\otimes
\mathcal{B}(\ell_2))}.
\end{eqnarray*}
Then we use the weak type $(1,1)$ boundedness of triangular
truncations to get
\begin{eqnarray*}
\mathsf{A}_\beta & \lesssim & \Big\| \sum_{k=1}^\infty \Big(
\sum_{j=1}^\infty \xi_{jk} d\beta_j \Big) \otimes e_{1k}
\Big\|_{L_1(\M \bar\otimes \mathcal{B}(\ell_2))} \\
& \le & \Big( \sup_{k \ge 1} \sum_{m=1}^\infty |\xi_{km}|^2
\Big)^\frac12 \, \sum_{j=1}^\infty \|d\beta_j\|_1 \lesssim \,
\sup_{n \ge 1} \|f_n\|_1.
\end{eqnarray*}
The last inequality follows from our hypothesis and Gundy's
decomposition. \fin

\subsection{$\mathrm{BMO}$ estimate, interpolation and duality}

In this paragraph we prove Theorem A2. The key for the BMO
estimate is certain commutation relation which can not be
exploited in $L_p$ for finite $p$. Namely, we have
$$\sum_{m=1}^\infty T_mf \otimes e_{1m} = \sum_{k=1}^\infty
\Big( \underbrace{\sum_{m=1}^\infty \xi_{km} \mathbf{1}_\M \otimes
e_{1m}}_{\xi_k^r} \Big) \big( df_k \otimes
\mathbf{1}_{\mathcal{B}(\ell_2)} \big) = \sum_{k=1}^\infty d \big(
\xi_k^r (f \otimes \mathbf{1}_{\mathcal{B}(\ell_2)}) \big)_k$$
where the last martingale difference is considered with respect to
the filtration $\mathcal{R}_n = \M_n \bar\otimes
\mathcal{B}(\ell_2)$ of $\mathcal{R}$. Note that $\xi_k^r$
commutes with $df_k \otimes \mathbf{1}_{\mathcal{B}(\ell_2)}$ and
therefore we find that
\begin{eqnarray*}
\Big\| \sum_{m=1}^\infty T_mf \otimes e_{1m}
\Big\|_{\mathrm{BMO}_{\mathcal{R}}^r}^2 & = & \Big\| \sup_{n \ge
1} \sum_{k \ge n} \mathsf{E}_n \big( (df_k \otimes
\mathbf{1}_{\mathcal{B}(\ell_2)}) \xi_k^r {\xi_k^r}^*(df_k \otimes
\mathbf{1}_{\mathcal{B}(\ell_2)})^* \big) \Big\|_{\mathcal{R}} \\
& \le & \Big( \sup_{k \ge 1} \sum_{m=1}^{\infty} |\xi_{km}|^2
\Big) \big\| f \otimes \mathbf{1}_{\mathcal{B}(\ell_2)}
\big\|_{\mathrm{BMO}_\mathcal{R}^r}^2 \\ [7pt] & \lesssim & \big\|
f \otimes \mathbf{1}_{\mathcal{B}(\ell_2)}
\big\|_{\mathrm{BMO}_\mathcal{R}^r}^2 \ \le \
\|f\|_{\mathrm{BMO}_\M^r}^2 \ \le \ \|f\|_\infty^2.
\end{eqnarray*}
Similarly, $\| \sum_m T_mf \otimes
e_{1m}\|_{\mathrm{BMO}_{\mathcal{R}}^c} \lesssim
\|f\|_{\mathrm{BMO}_\M^c} \le \|f\|_\infty$ so that $$\Big\|
\sum_{m=1}^\infty T_mf \otimes e_{1m}
\Big\|_{\mathrm{BMO}(\mathcal{R})} \lesssim
\|f\|_{\mathrm{BMO}(\M)} \le \sup_{n \ge 1} \|f_n\|_\infty.$$ The
estimate for $\sum_m T_mf \otimes e_{m1}$ is entirely analogous.
This gives the $L_\infty - \mathrm{BMO}$ estimate, or even better
the $\mathrm{BMO}-\mathrm{BMO}$ one. Note that we make crucial use
of the identity $\|f \otimes
\mathbf{1}_{\mathcal{B}(\ell_2)}\|_\infty = \|f\|_\infty$!
However, $\|f \otimes
\mathbf{1}_{\mathcal{B}(\ell_2)}\|_{L_p({\mathcal R})} \neq
\|f\|_{L_p({\mathcal M})}$ for $p$ finite. That is why we can not
reduce the $L_p$ estimate to the commutative case.

With the weak type $(1,1)$ and the ${\mathrm {BMO}}$ estimates in
hand, we may follow by interpolation. Namely, since the case $p=2$
is trivial, we interpolate for $1 < p < 2$ following
Randrianantoanina's argument \cite{R3} and for $2 < p < \infty$
following Junge and Musat \cite{JMu,Mu}. This gives rise to
$$\Big\| \sum_{m=1}^\infty T_mf \otimes \delta_m \Big\|_{L_p({\mathcal
M}; \ell_{rc}^2)} \le c_p \, \|f\|_p$$ for all $1<p<\infty$.
Assuming further that $\sum_m |\xi_{km}|^2 = \gamma_k \sim 1$, we
find
\begin{eqnarray*}
\|f\|_p & = & \sup_{\|g\|_{p'} \le 1} \, \sum_{k=1}^\infty
\big\langle df_k, dg_k \big\rangle \\ & = & \sup_{\|g\|_{p'} \le
1} \, \sum_{k=1}^\infty \frac{1}{\gamma_k} \sum_{m=1}^\infty
\big\langle \xi_{km} df_k, \xi_{km} dg_k \big\rangle \\
& = & \sup_{\|g\|{p'} \le 1} \, \sum_{m=1}^\infty \Big\langle
\sum_{k=1}^\infty \xi_{km}
df_k, \sum_{k=1}^\infty \frac{\xi_{km}}{\gamma_k} dg_k \Big\rangle \\
& \le & \Big\| \sum_{m=1}^\infty T_mf \otimes \delta_m
\Big\|_{L_p({\mathcal M};\ell_{rc}^2)} \sup_{\|g\|_{p'} \le 1} \,
\Big\| \sum_{m=1}^\infty \sum_{k=1}^\infty
\frac{\xi_{km}}{\gamma_k} dg_k \otimes \delta_m
\Big\|_{L_{p'}({\mathcal M};\ell_{rc}^2)}.
\end{eqnarray*}
Therefore, the reverse inequality follows by duality whenever
$\sum_m |\xi_{km}|^2 \sim 1$. \fin

\begin{remark}
\emph{The aim of Theorem A2 is $$\Big\| \Big( \sum_{m=1}^\infty
(T_m \hskip-1pt f)(T_m \hskip-1pt f)^* \Big)^\frac12 \Big\|_p +
\Big\| \Big( \sum_{m=1}^\infty (T_m \hskip-1pt f)^*(T_m \hskip-1pt
f) \Big)^\frac12 \Big\|_p \lesssim \, c_p \, \|f\|_p$$ for $2 < p
< \infty$, since the remaining inequalities follow from it and
Theorem A1. Nevertheless, a direct argument (not including the BMO
estimate) is also available for $p \in 2 \hskip1pt \mathbb{Z}$.
Indeed, by Khintchine and Burkholder-Gundy inequalities
\begin{eqnarray*}
\Big\| \sum_{m=1}^\infty T_mf \otimes \delta_m
\Big\|_{L_p({\mathcal M};\ell_{rc}^2)}^p & \sim & \int_\Omega
\Big\| \sum_{k=1}^\infty \Big( \underbrace{\sum_{m=1}^\infty
\xi_{km} r_m(w)}_{\xi_k(w)} \Big) df_k \Big\|_{L_p(\M)}^p \,
d\mu(w) \\ & \sim & \int_\Omega \Big\| \sum_{k=1}^\infty \xi_k(w)
df_k \otimes \delta_k \Big\|_{L_p(\M; \ell_{rc}^2)}^p \, d\mu(w).
\end{eqnarray*}
If $p=2j$, we use H\"older and again Khintchine + Burkholder-Gundy
to obtain
\begin{eqnarray*}
\lefteqn{\hskip-20pt \tau \int_\Omega \Big(
\sum_{k=1}^\infty |\xi_k|^2 |df_k|^2 \Big)^\frac{p}{2} \, d\mu}
 \\ & = & \sum_{k_1, k_2, \ldots, k_j=1}^\infty
\Big( \int_\Omega \prod_{s=1}^j |\xi_{k_s}|^2 \, d\mu \Big) \,
\tau \Big[ |df_{k_1}|^2 |df_{k_2}|^2 \cdots |df_{k_j}|^2 \Big] \\
& \le & \sum_{k_1, k_2, \ldots, k_j=1}^\infty \prod_{s=1}^j \Big(
\int_\Omega |\xi_{k_s}|^{2j} \, d\mu \Big)^{\frac1j} \tau \Big[
|df_{k_1}|^2 |df_{k_2}|^2 \cdots |df_{k_j}|^2\Big]\\ & \lesssim &
\Big( \sup_{k \ge 1} \sum_{m=1}^\infty |\xi_{km}|^2
\Big)^{\frac{1}{2}} \tau \Big( \sum_{k=1}^\infty |df_k|^2
\Big)^{\frac{p}{2}} \ \le \ c_p^p \, \|f\|_p^p.
\end{eqnarray*}
The row term is estimated in the same way. This completes the
argument. \fin}
\end{remark}

\begin{remark}
\emph{Our results so far and the implications for semigroups
explored in the next paragraph can be regarded as an alternative
argument in producing Littlewood-Paley inequalities from
martingale inequalities. Namely, the key result is Gundy's
decomposition while in \cite[Chapter IV]{St2} the main ideas are
based on Stein's inequality for martingales, which is not
necessary from our viewpoint.}
\end{remark}

\begin{remark}
\emph{The adjoint of $\sum_m T_m \otimes \delta_m$ is
$$\widehat{g} = \sum_{m=1}^\infty g^m \otimes \delta_m \in
L_p(\M; \ell_{rc}^2) \mapsto \sum_{m=1}^\infty T_m g^m \in
L_p(\M).$$ Letting $\widehat{\xi}_k = \sum_m \xi_{km} \otimes
\delta_m$, this mapping can be formally written as
$$\sum_{m=1}^\infty T_m g^m = \sum_{k=1}^\infty \sum_{m=1}^\infty
\xi_{km} dg_k^m = \sum_{k=1}^\infty \big\langle \widehat{\xi}_k, d
\widehat{g}_k \big\rangle.$$ In other words, we find a martingale
transform with row/column valued coefficients and martingale
differences. In this case, we have again noncommuting coefficients
and Theorem A2 gives the right $L_p$ estimate.}
\end{remark}

\subsection{Rota's dilation theorem and pseudo-localization}

We show that the key new ingredient to produce weak type
inequalities for square functions associated to a large family of
semigroups is a certain pseudo-localization estimate. This applies
in particular for semicommutative Calder\'on-Zygmund operators and
illustrates the difference between Theorems A1 and B1. Our
argument uses Rota's theorem \cite{Ro} in conjunction with Stein's
ergodic averages \cite{St2}, we find it quite transparent.

According to \cite{JLX}, we will say that a bounded operator $T:
\M \to \M$ satisfies \emph{Rota's dilation property} if there
exists a von Neumann algebra $\mathcal{R}$ equipped with a
normalized trace, a normal unital faithful $*$-representation
$\pi: \M \to \mathcal{R}$ which preserves trace, and a decreasing
sequence $(\mathcal{R}_m)_{m \ge 1}$ of von Neumann subalgebras of
$\mathcal{R}$ such that $T^m = \mathbb{E} \circ \mathcal{E}_m
\circ \pi$ for any $m \ge 1$ and where $\mathcal{E}_m: \mathcal{R}
\to \mathcal{R}_m$ is the canonical conditional expectation and
$\mathbb{E}: \mathcal{R} \to \M$ is the conditional expectation
associated with $\pi$. By Rota's theorem, $T^2$ satisfies it
whenever $\M$ is commutative and $T: \M \to \M$ is a normal unital
positive self-adjoint operator. This was used by Stein \cite{St2}
for averages of discretized symmetric diffusion commutative
semigroups and also applies in the semicommutative setting of
\cite{Me}. We also know from \cite{JLX} that the noncommutative
Poisson semigroup on the free group fits in.

The problem that we want to study is the following. Assume that
$T: \M \to \M$ is a normal unital completely positive self-adjoint
operator satisfying Rota's dilation property. Given $f \in
L_1(\M)$ and $m \ge 0$, set $$\Sigma_m \hskip-1pt f =
\frac{1}{m+1} \sum_{k=0}^m T^k \! f \quad \mbox{and} \quad
\Gamma_m \hskip-1pt f = \Sigma_m \hskip-1pt f - \Sigma_{m-1}
\hskip-1ptf.$$ What can we say about the inequality $$\Big\|
\sum_{m=1}^\infty \sqrt{m} \ \Gamma_m \hskip-1pt f \otimes
\delta_m \Big\|_{L_{1,\infty}(\M;\ell_{rc}^2)} \lesssim \|f\|_1?$$
Of course, the norm in $L_{1,\infty}(\M;\ell_{rc}^2)$ denotes the
one used in Theorem A1.

This might be related to weak type estimates for general symmetric
diffusion semigroups satisfying Rota's dilation property. It seems
though that the classical method \cite{St2} only works for $p>1$,
see \cite{JLX,JX2} for noncommutative forms of Stein's fractional
averages. Nevertheless there are concrete cases, like the
noncommutative Poisson semigroup on the free group, where more
information is available. The idea is to prove first its
martingale analog and apply then Rota's property to recognize
pseudo-localization as the key new ingredient. The martingale
inequality below is the weak type (1,1) extension of
\cite[Proposition 10.8]{JLX}.

\begin{corollary}
Let $(\M_n)_{n \ge 1}$ be either an increasing or decreasing
filtration of the von Neumann algebra $\M$. Given $f \in L_1(\M)$,
we set $f_n = \mathcal{E}_n(f)$ and define the following operators
associated to $f$ for $m \ge 0$ $$\sigma_m \hskip-1pt f =
\frac{1}{m+1} \sum_{k=0}^m f_k \quad \mbox{and} \quad \gamma_m
\hskip-1pt f = \sigma_m \hskip-1pt f - \sigma_{m-1} \hskip-1ptf.$$
Then, there exists a decomposition $\sqrt{m} \, \gamma_m
\hskip-1pt f = \alpha_m \hskip-1pt f + \beta_m \hskip-1pt f$ such
that $$\Big\| \Big( \sum_{m=1}^\infty (\alpha_m \hskip-1pt f)
(\alpha_m \hskip-1pt f)^* \Big)^\frac12 \Big\|_{1,\infty} + \Big\|
\Big( \sum_{m=1}^\infty (\alpha_m \hskip-1pt f)^* (\alpha_m
\hskip-1pt f) \Big)^\frac12 \Big\|_{1,\infty} \lesssim \|f\|_1.$$
\end{corollary}

\dem Theorem A1 still holds for decreasing filtrations $(\M_n)_{n
\ge 1}$. Indeed, it suffices to observe that any finite (reverse)
martingale $(f_1, f_2, \ldots, f_n, f_n, f_n, \ldots)$ may be
rewritten as a finite ordinary martingale by reversing the order,
and that this procedure does not affect the arguments in the proof
of Theorem A1, details are left to the reader. On the other hand,
as in \cite{JLX,St2} it is easily checked that $$\sqrt{m} \,
\gamma_m \hskip-1pt f = \sum_{k=1}^m \frac{k}{\sqrt{m} (m+1)} \,
df_k.$$ Therefore, the result follows from Theorem A1 since
$\xi_{km} = \delta_{k \le m} \frac{k}{\sqrt{m}(m+1)}$. \fin

Now we are in a position to study the norm of $\sum_m \sqrt{m} \
\Gamma_m \hskip-1pt f \otimes \delta_m$ in
$L_{1,\infty}(\M;\ell_{rc}^2)$. Namely, since we have assumed that
$T$ satisfies Rota's dilation property, we have $\sqrt{m} \
\Gamma_m f = \sqrt{m} \ \mathbb{E} \circ \gamma_m \circ \pi(f)$.
The presence of $\pi$ is harmless since it just means that we are
working in a bigger algebra. Moreover, when working in $L_p$ the
conditional expectation $\mathbb{E}$ is a contraction so that we
may eliminate it. However, this is not the case in $L_{1,\infty}$
and we have to review the arguments in the proof of Theorem A1 for
the operator
$$T_m \hskip-1pt f = \sqrt{m} \ \mathbb{E} \circ \gamma_m
\hskip-1pt f = \sum_{k=1}^m \frac{k}{\sqrt{m}(m+1)} \, \mathbb{E}
\, df_k = \sum_{k=1}^\infty \xi_{km} \, \mathbb{E} \, df_k.$$ We
consider the decomposition $T_m \hskip-1pt f = A_m \hskip-1pt f +
B_m \hskip-1pt f$ with $$A_m \hskip-1pt f = \sum_{k=1}^\infty
\xi_{km} \Delta_r (\mathbb{E} \, df_k) \quad \mbox{and} \quad B_m
\hskip-1pt f = \sum_{k=1}^\infty \xi_{km} \Delta_c (\mathbb{E} \,
df_k).$$ Following the proof of Theorem A1, we are reduced to
estimate the three square functions associated to the terms
$\alpha, \beta$ and $\gamma$ in Gundy's decomposition. The
$\alpha$ and $\beta$ terms are estimated in the same way, since
the weak $L_1$ and the $L_2$ boundedness of $\Delta_{r\ell}
\otimes id_{\mathcal{B}(\ell_2)}$ is also satisfied by
$\Delta_{r\ell} \mathbb{E} \otimes id_{\mathcal{B}(\ell_2)}$.

\noindent \textbf{Conclusion.} The only term that can not be
estimated from the argument in Theorem A1 is the $\gamma$-term.
Nevertheless, the key to estimate that term is the fact that it is
supported by a sufficiently small projection. More specifically,
it is easily checked that we have $\mathrm{supp}^* df_k \le
\mathbf{1}_\M - w_\ell$ for $k \ge 1$ and $\lambda \tau
(\mathbf{1}_\M - w_\ell) \le 2 \|f\|_1$. According to \cite[Remark
5.1]{P2}, this gives $$df_k = (\mathbf{1}_\M - w_\ell)df_k +
df_k(\mathbf{1}_\M - w_\ell) - (\mathbf{1}_\M -
w_\ell)df_k(\mathbf{1}_\M - w_\ell).$$ Moreover, since $w_\ell
\pi_k = \pi_k w_\ell = \pi_k$ for $k \le \ell$, we find
\begin{eqnarray*}
\sum_{k=1}^\infty \xi_{km} \Delta_{r\ell}(\mathbb{E} df_k) & = &
\Delta_{r\ell} \Big[ w_\ell \mathbb{E} \Big( w_\ell^\perp
\Lambda_m f + \Lambda_m f w_\ell^\perp - w_\ell^\perp \Lambda_m f
w_\ell^\perp \Big) w_\ell \Big],
\end{eqnarray*}
where $\Lambda_m f = \sum_k \xi_{km} df_k$ and $w_\ell^\perp =
\mathbf{1}_\M - w_\ell$. That is why we obtain a zero term in the
martingale case, with $\mathbb{E}$ being the identity map. In the
general case, this indicates that we shall be able to obtain
satisfactory inequalities whenever $\mathbb{E}$ \emph{almost}
behaves as a local map, in the sense that \emph{respects} the
supports. As we shall see in the next section, when dealing with
operator-valued Calder\'on-Zygmund operators, the role of $\gamma$
will be played by the off-diagonal terms of the good and bad parts
of Calder\'on-Zygmund decomposition. The pseudo-localization
estimate needed for the bad part is standard, while the one for
the good part requires some almost-orthogonality methods described
in Appendix A.

As far as we know, this pseudo-localization property is unknown
for all the \emph{fully} noncommutative Calder\'on-Zygmund-type
operators in the literature. Particularly it would be very
interesting to know the behavior of the noncommutative Poisson
semigroup on the free group. This leads us to formulate it as a
problem for the interested reader.

\noindent \textbf{Problem.} Let us consider the Poisson semigroup
$$T_t(\lambda(g)) = e^{-t|g|} \lambda(g)$$ given by the length
function in the free group von Neumann algebra. It is known
\cite{JLX,JX2} that square and maximal functions are bounded maps
on $L_p$. Are these mappings of type $(1,1)$?

\subsection{Further comments}

We have seen so far several particular cases of Theorems A1 and
A2. Namely, in the Introduction we mentioned noncommutative
martingale transforms and square functions as well as
\emph{shuffled} square functions where the martingale differences
are grouped according to an arbitrary partition of $\mathbb{N}$.
In the previous paragraph we have seen that Stein's ergodic
averages, in connection with semigroups satisfying Rota's
property, also fall in the framework of Theorems A1 and A2. In
this paragraph, we indicate further applications and
generalizations of Theorems A1 and A2.

\noindent A. \emph{Multi-indexed coefficients/martingales.}
Replacing $df_k$ in our main results by Rademacher variables or
free generators clearly gives rise to new Khintchine type
inequalities. The iteration of Khintchine inequalities was the
basis in \cite{PP} for some multilinear generalizations that where
further explored in \cite{JPX,RX}. Although a detailed analysis of
these methods in the context of our results is out of the scope of
this paper, let us mention two immediate consequences.

In the iteration of Khintchine type inequalities, it is sometime
quite interesting to being able to dominate the cross terms that
appear by the row/column terms and thereby reduce it to a standard
Khintchine type inequality. Our result is \emph{flexible} enough
to produce such estimates.

\begin{corollary} \label{Cross}
Assume that $$\sum_{m=1}^\infty |\rho_{km}|^2 \sim 1 \sim
\sum_{n=1}^\infty |\eta_{kn}|^2\quad \mbox{for all} \quad k \ge
1.$$ Let $T_{mn}f = \sum_k \rho_{km} \eta_{kn} df_k$. Then, if $2
< p < \infty$, we have $$\Big\| \sum_{m,n=1}^\infty T_{mn}f
\otimes e_{m,n} \Big\|_p \lesssim \Big\| \sum_{m,n=1}^\infty
T_{mn}f \otimes e_{1,mn} \Big\|_p + \Big\| \sum_{m,n=1}^\infty
T_{mn}f \otimes e_{mn,1} \Big\|_p.$$ We may also replace $e_{m,n}$
by $e_{n,m}$. If $1 < p < 2$, certain dual inequalities hold.
\end{corollary}

\dem Let us consider the spaces
\begin{eqnarray*}
\mathcal{K}_p^1 & = & R_p \otimes_h R_p + C_p \otimes_h C_p, \\
\mathcal{J}_p^1 & = & R_p \otimes_h R_p \hskip0.5pt \cap
\hskip0.5pt C_p \otimes_h C_p, \\
\mathcal{K}_p^2 & = & R_p \otimes_h R_p + R_p \otimes_h C_p + C_p
\otimes_h R_p + C_p \otimes_h C_p, \\
\mathcal{J}_p^2 & = & R_p \otimes_h R_p \hskip0.5pt \cap
\hskip0.5pt R_p \otimes_h C_p \hskip0.5pt \cap \hskip0.5pt C_p
\otimes_h R_p \hskip0.5pt \cap \hskip0.5pt C_p \otimes_h C_p.
\end{eqnarray*}
The assertion for $2 < p < \infty$ gives that the norm of
$(T_{mn}f)$ in $L_p(\M; \mathcal{J}_p^2)$ is bounded above by the
norm in $L_p(\M; \mathcal{J}_p^1)$. The dual formulation just
means that $$\big\| (T_{mn}f)_{m,n \ge 1}
\big\|_{L_p(\M;\mathcal{K}_p^1)} \lesssim \big\| (T_{mn}f)_{m,n
\ge 1} \big\|_{L_p(\M;\mathcal{K}_p^2)}.$$ The proofs of both inequalities are
similar, thus we can assume $2 < p < \infty$. Since
$$\sum_{m,n=1}^\infty |\rho_{km} \eta_{kn}|^2 \sim 1,$$ we may
apply Theorem A2 for $\xi_{k,(m,n)} = \rho_{km} \eta_{kn}$ and
obtain that the norm of $(T_{mn}f)$ in $L_p(\M;\mathcal{J}_p^1)$
is equivalent to the norm of $f$ in $L_p(\M)$. We claim that the
same equivalence holds for $L_p(\M; \mathcal{J}_p^2)$. Indeed,
applying Theorem A2 first in the variable $m$ and then in $n$, we
obtain
\begin{eqnarray*}
\lefteqn{\hskip-25pt \Big\| \sum_{m,n=1}^\infty T_{mn}f \otimes
e_{m,n} \Big\|_p} \\ & = & \Big\| \sum_{m=1}^\infty \Big[
\sum_{k=1}^\infty \rho_{km} \Big( \sum_{n=1}^\infty \eta_{kn}
e_{1,n} \Big) \otimes df_k \Big] \otimes e_{m,1} \Big\|_p \\ &
\lesssim & \Big\| \sum_{n=1}^\infty \Big( \sum_{k=1}^\infty
\eta_{kn} df_k \Big) \otimes e_{1,n} \Big\|_p \ \lesssim \
\|f\|_p.
\end{eqnarray*}
The same argument for $e_{n,m}$ instead of $e_{m,n}$ applies and
the result follows. \fin

Note that Corollary \ref{Cross} also applies for $d$ indices. In
this case all the terms are dominated by the row/column terms
$e_{1,m_1 \cdots m_d}$ and $e_{m_1 \cdots m_d,1}$. Another
multilinear form of Theorem A2 is given by considering
multi-indexed martingales. We refer to \cite[Section 4.2]{P1} for
the definition of a multi-indexed martingale. Like it was pointed
in \cite{P1}, successive iterations of Theorem A2 give rise to a
generalization of it for multi-indexed martingales.

\noindent B. \emph{Commuting operator coefficients.} A natural
question is whether Theorems A1 and A2 still hold when the
coefficients $\xi_{km}$ are operators instead of scalars. In view
of \cite{R1}, we must impose certain commuting condition of the
coefficients $\xi_{km}$, we refer to \cite[Section 6]{P2} for more
details on this topic. We state below the precise statement for
operator coefficients. It is not difficult to check that the same
arguments apply to this more general case, details are left to the
reader.

\begin{corollary}
Assume that $\xi_{km} \in \M_{k-1} \cap \M'$ and $$\sup_{k \ge 1}
\ \Big\| \sum_{m=1}^\infty \xi_{km} \xi_{km}^* \Big\|_\M + \Big\|
\sum_{m=1}^\infty \xi_{km}^* \xi_{km} \Big\|_\M \lesssim 1,$$
where $\M'$ denotes the commutant of $\M$. Then, Theorems A1 and
A2 still hold.
\end{corollary}

\begin{remark}
\emph{It is finally worth mentioning that Junge and K\"ostler have
recently developed an $H_p$ theory of noncommutative martingales
for continuous filtrations \cite{JK}. Our results could also be
studied in such setting, although it is again out of the scope of
the paper.}
\end{remark}

\section{Calder\'on-Zygmund operators} \label{S2}

In this section, we prove Theorems B1 and B2. We recall the
definition of the von Neumann algebra $(\Mn,\varphi)$ from the
Introduction. Note that $L_p(\Mn)$ is the Bochner space of $L_p$
functions on $\R^n$ with values in $L_p(\M)$. We shall need some
additional terminology. The size of a cube $Q$ in $\R^n$ is the
length $\ell(Q)$ of  its edges. Given $k \in \Z$, we use $\Q_k$
for the set of dyadic cubes of size $1/2^{k}$. If $Q$ is a dyadic
cube and $f: \R^n \to L_p(\M)$, we set $$f_Q = \frac{1}{|Q|}
\int_Q f(y) \, dy.$$ Let $(\mathsf{E}_k)_{k \in \Z}$ denote this
time the family of conditional expectations associated to the
classical dyadic filtration on $\R^n$. $\mathsf{E}_k$ will also
stand for the tensor product $\mathsf{E}_k \otimes id_\M$ acting
on $\Mn$. If $1 \le p < \infty$ and $f \in L_p(\Mn)$,
$$\mathsf{E}_k(f) = f_k = \sum_{Q \in \Q_k}^{\null} f_Q 1_Q.$$
$(\Mn_k)_{k \in \Z}$ denotes the corresponding filtration $\Mn_k =
\mathsf{E}_k(\Mn)$. We use $\widehat{Q}$ for the dyadic father of
a dyadic cube $Q$, the dyadic cube containing $Q$ with double
size. Given $\delta > 1$, the $\delta$-concentric father $\delta
\hskip1pt Q$ of $Q$ is the cube concentric with $Q$ satisfying
$\ell(\delta \hskip1pt Q) = \delta \hskip1pt \ell(Q)$. Given $f:
\R^n \to \C$, let $df_k$ denote the $k$-th martingale difference
with respect to the dyadic filtration. That is, $$df_k = \sum_{Q
\in \Q_k} \big( f_Q - f_{\widehat{Q}} \big) 1_Q.$$ Let
$\mathcal{R}_k$ be the class of sets in $\R^n$ being the union of
a family of cubes in $\Q_k$. Given such an $\mathcal{R}_k$-set
$\Omega = \bigcup_j Q_j$, we shall work with the dilations $9
\Omega = \bigcup_j 9 Q_j$, where $9 Q$ denotes the $9$-concentric
father of $Q$.

The right substitute of Gundy's decomposition in our new setting
will be the noncommutative form of Calder\'on-Zygmund
decomposition. Again, the main tool is Cuculescu's construction
associated (this time) to the filtration $(\Mn_k)_{k \in \Z}$. Let
us consider the dense subspace $$\Mn_{c,+} = L_1(\Mn) \cap \Big\{
f: \R^n \to \M \, \big| \ f \in \Mn_+, \
\overrightarrow{\mathrm{supp}} \hskip1pt f \ \ \mathrm{is \
compact} \Big\} \subset L_1(\Mn).$$ Here
$\overrightarrow{\mathrm{supp}}$ means the support of $f$ as a
vector-valued function in $\R^n$. In other words, we have
$\overrightarrow{\mathrm{supp}} \hskip1pt f = \mathrm{supp}
\hskip1pt \|f\|_\M$. We employ this terminology to distinguish
from $\mathrm{supp} \, f$, the support of $f$ as an operator in
$\Mn$. Any function $f \in \Mn_{c,+}$ gives rise to a martingale
$(f_n)_{n \in \Z}$ with respect to the dyadic filtration and we
may consider the Cuculescu's sequence $(q_k(\lambda))_{k \in \Z}$
associated to $(f,\lambda)$ for any $\lambda >0$. Since $\lambda$
will be fixed most of the time, we will shorten the notation by
$q_k$ and only write $q_k(\lambda)$ when needed. Define the
sequence $(p_k)_{k \in \Z}$ of disjoint projections $$p_k =
q_{k-1}-q_k.$$ As noted in \cite{P2}, we have $q_k =
\mathbf{1}_\Mn$ for $k$ small enough and $$\sum_{k \in \Z} p_k =
\mathbf{1}_\Mn - q \quad \mbox{with} \quad q = \bigwedge_{k \in
\Z} q_k.$$

\begin{CZDecomposition}
Let $f \in \Mn_{c,+}$ and let $\lambda$ be a positive number.
Then, $f$ can be decomposed $f = g_d + g_\mathit{off} + b_d +
b_\mathit{off}$ as the sum of four functions
$$\begin{array}{rclcrcl} \displaystyle g_d & = & \displaystyle
qfq + \sum_{k \in \Z} p_k f_k p_k, & \quad & g_\mathit{off} & = &
\displaystyle \sum_{i \neq j} p_i f_{i \vee j} p_j \ + \ q f
q^\perp + q^\perp f q, \\ [15pt] b_d & = & \displaystyle \sum_{k
\in \Z} p_k \hskip1pt (f - f_k) \hskip1pt p_k, & \quad &
b_{\mathit{off}} & = & \displaystyle \sum_{i \neq j} p_i (f-f_{i
\vee j}) p_j,
\end{array}$$ where $i \vee j = \max (i,j)$ and $q^\perp =
\mathbf{1}_\Mn - q$. Moreover, we have the diagonal estimates
$$\Big\| qfq + \sum_{k \in \Z} p_k f_k p_k \Big\|_2^2 \le 2^n
\lambda \, \|f\|_1 \quad \mbox{and} \quad \sum_{k \in \Z} \big\|
p_k (f-f_k) p_k \big\|_1 \le 2 \, \|f\|_1.$$
\end{CZDecomposition}

\begin{remark}
\emph{Some comments:}
\begin{itemize}
\item \emph{Let $m_\lambda$ be the larger integer with
$q_{m_\lambda}(\lambda) = \mathbf{1}_\Mn$.}
\item \emph{There exist weaker off-diagonal estimates for $g$ and
$b$, see \cite[Appendix B]{P2}.}
\end{itemize}
\end{remark}

\begin{remark} \label{goff}
\emph{We have $$g_{\mathit{off}} = \sum_{s=1}^\infty
\sum_{k=m_\lambda+1}^\infty p_k df_{k+s} q_{k+s-1} + q_{k+s-1}
df_{k+s} p_k = \sum_{s=1}^\infty \sum_{k=m_\lambda+1}^\infty
g_{k,s} = \sum_{s=1}^\infty g_{(s)}.$$ Moreover, it is easily
checked that
$$\sup_{s \ge 1} \|g_{(s)}\|_2^2 = \sup_{s \ge 1}
\sum_{k=m_\lambda+1}^\infty \|g_{k,s}\|_2^2 \lesssim \lambda \,
\|f\|_1$$ and that $\mathrm{supp}^* {dg_{(s)}}_{k+s} =
\mathrm{supp}^* g_{k,s} \le p_k \le \mathbf{1}_\Mn-q_k$, see
\cite{P2} for further details.}
\end{remark}

We need one more preliminary result. Given $\lambda > 0$, we adopt
the terminology from \cite{P2} and write $q_k(\lambda) = \sum_{Q
\in \Q_k} \xi_Q 1_Q$ with $\xi_Q$ projections in $\M$. Thus, since
we are assuming that $q_{m_\lambda}(\lambda)=\mathbf{1}_\Mn$, we
may write $$\mathbf{1}_\Mn - q_k(\lambda) = \sum_{s=m_\lambda+1}^k
\sum_{Q \in \Q_s} (\xi_{\widehat{Q}} - \xi_Q) 1_Q.$$ An
$\R^n$-dilated version (by a factor $9$) of it is given by
$$\mathrm{supp} \, \psi_k(\lambda) \quad \mbox{with} \quad
\psi_k(\lambda) = \sum_{s=m_\lambda+1}^k \sum_{Q \in \Q_s}
(\xi_{\widehat{Q}} - \xi_Q) 1_{9Q},$$ the support projection of
$\psi_k(\lambda)$. The result below is proved in \cite[Lemma
4.2]{P2}.

\begin{lemma} \label{keylem}
Let us set $$\zeta(\lambda) = \bigwedge_{k \in \Z}
\zeta_k(\lambda) \quad \mbox{with} \quad \zeta_k(\lambda) =
\mathbf{1}_\Mn - \mathrm{supp} \, \psi_k(\lambda).$$ Then,
$\zeta(\lambda)$ is a projection in $\Mn$ and we have
\begin{itemize}
\item[i)] $\lambda \hskip1pt \varphi \big( \1_\Mn - \zeta(\lambda)
\big) \le 9^n \hskip1pt \|f\|_1$.
\item[ii)] If $Q_0$ is any dyadic cube, then we have
$$\zeta(\lambda) \big( \mathbf{1}_\M \otimes 1_{9Q_0} \big) \le
\big( \1_\M - \xi_{\widehat{Q}_0} + \xi_{Q_0} \big) \otimes
1_{9Q_0}.$$ In particular, it can be deduced that $\zeta(\lambda)
\big( \mathbf{1}_\M \otimes 1_{9Q_0} \big) \le \xi_{Q_0} \otimes
1_{9Q_0}$.
\end{itemize}
\end{lemma}

\subsection{A pseudo-localization result}
\label{2.1}

Now we show how to transfer the result in Appendix A to the
noncommutative setting. We shall need the weak notion
$\mathrm{supp}^*$ of support projection introduced before the
statement of Gundy's decomposition above. It is easily seen that
$\mathrm{supp}^* f$ is the smallest projection $p$ in $\Mn$
satisfying the identity $f = p f + f p - p f p$. We shall prove
the following result.

\begin{theorem} \label{pseudolocal}
Given a Hilbert space $\mathcal{H}$, consider a kernel $k: \R^{2n}
\setminus \Delta \to \mathcal{H}$ which satisfies the
size/smoothness conditions imposed in the Introduction and
formally defines the Calder\'on-Zygmund operator $$T \! f(x) =
\int_{\R^n} k(x,y) f(y) \, dy.$$ Assume further that $T: L_2(\Mn)
\to L_2(\Mn;\mathcal{H}_{oh})$ is of norm $1$ and fix a positive
integer $s$. Given $f \in L_2(\Mn)$ and $k \in \Z$, let us
consider any projection $q_k$ in $\Mn_k$ satisfying that $\1_\Mn -
q_k$ contains $\mathrm{supp}^* df_{k+s}$ as a subprojection. If we
write $$q_k = \sum_{Q \in \Q_k} \xi_Q 1_Q$$ with $\xi_Q$
projections, we may further consider the projection
$$\zeta_{f,s} = \bigwedge_{k \in \Z} \Big( \1_\Mn -
\bigvee_{Q \in \Q_k} (\1_\M - \xi_Q) \hskip1pt 1_{9Q} \Big).$$
Then we have the following localization estimate in
$L_2(\Mn;\mathcal{H}_{oh})$
$$\Big( \int_{\R^n} \left\| \big[ \zeta_{f,s} \, T \! f
\, \zeta_{f,s} \big] (x) \right\|_{L_2(\M;\mathcal{H}_{oh})}^2 \,
dx \Big)^{\frac12} \le \mathrm{c}_{n,\gamma} s \hskip1pt 2^{-
\gamma s/2} \Big( \int_{\R^n} \tau \big( |f(x)|^2 \big) \, dx
\Big)^{\frac12}.$$
\end{theorem}

\dem We shall reduce this result to its commutative counterpart in
Appendix A below. \hskip3pt According to the shift condition
imposed \hskip3pt $\mathrm{supp}^* df_{k+s} \prec \1_\Mn - q_k$,
we have
$$df_{k+s} = q_k^\perp df_{k+s} + df_{k+s} q_k^\perp - q_k^\perp
df_{k+s} q_k^\perp$$ where we write $q_k^\perp = \1_\Mn - q_k$ for
convenience. On the other hand, let
\begin{eqnarray}
\zeta_k = \1_\Mn -
\bigvee_{Q \in \Q_k} (\1_\M - \xi_Q) \hskip1pt 1_{9Q}, \label{zeta}
\end{eqnarray}
 so that
$\zeta_{f,s} = \bigwedge_k \zeta_k$. As in Lemma \ref{keylem}, it
is easily seen that $\1_\Mn - \zeta_k$ represents the
$\R^n$-dilated projection associated to $\1_\Mn - q_k$ with a
factor $9$. Let $\mathcal{L}_a$ and $\mathcal{R}_a$ denote the
left and right multiplication maps by the operator $a$. Let also
$\mathcal{LR}_a$ stand for $\mathcal{L}_a + \mathcal{R}_a -
\mathcal{L}_a \mathcal{R}_a$ Then our considerations so far and
the fact that $\mathcal{L}_{\zeta_k}, \mathcal{R}_{\zeta_k}$ and
$\mathcal{LR}_{q_k^{\perp}}$ commute with $\mathsf{E}_j$ for $j
\ge k$ give
\begin{eqnarray*}
\lefteqn{\zeta_{f,s} \, T \! f \, \zeta_{f,s}} \\ & = &
\mathcal{L}_{\zeta_{f,s}} \mathcal{R}_{\zeta_{f,s}} \Big( \summ_k
\mathsf{E}_k T \Delta_{k+s} \mathcal{LR}_{q_k^\perp} + \summ_k (id
- \mathsf{E}_k) \mathcal{L}_{\zeta_k} \mathcal{R}_{\zeta_k} T
\mathcal{LR}_{q_k^\perp} \Delta_{k+s} \Big) (f).
\end{eqnarray*}
Now we claim that $$\mathcal{L}_{\zeta_k} \mathcal{R}_{\zeta_k} T
\mathcal{LR}_{q_k^\perp} = \mathcal{L}_{\zeta_k}
\mathcal{R}_{\zeta_k} T_{4 \cdot 2^{-k}}
\mathcal{LR}_{q_k^\perp}.$$ Indeed, this clearly reduces to see
$$\mathcal{L}_{\zeta_k} T
\mathcal{L}_{q_k^\perp} = \mathcal{L}_{\zeta_k} T_{4 \cdot 2^{-k}}
\mathcal{L}_{q_k^\perp} \qquad \mbox{and} \qquad
\mathcal{R}_{\zeta_k} T \mathcal{R}_{q_k^\perp} =
\mathcal{R}_{\zeta_k} T_{4 \cdot 2^{-k}}
\mathcal{R}_{q_k^\perp}.$$ By symmetry, we just prove the first
identity $$\mathcal{L}_{\zeta_k} T \mathcal{L}_{q_k^\perp} f (x)
\, = \sum_{Q \in \Q_k} \zeta_k(x) (\1_\M - \xi_Q) \int_Q k(x,y)
f(y) \, dy.$$ Assume that $x \in 9Q$ for some $Q \in \Q_k$, then
it follows as in Lemma \ref{keylem} (see the definition of the
projection $\zeta_k$) that $\zeta_k(x) \le \xi_Q$. In particular,
we deduce from the expression above that for each $y \in Q$ we
must have $x \in \R^n \setminus 9Q$. This implies $|x-y| \ge 4
\cdot 2^{-k}$ as desired. Finally, the operators $\mathcal{L},
\mathcal{R}$ and $\mathcal{LR}$ inside the bracket are clearly
absorbed by $f$ and $\zeta_{f,s}$. Thus, we obtain the identity
below
$$\zeta_{f,s} \, T \! f \, \zeta_{f,s} \, = \,
\mathcal{L}_{\zeta_{f,s}} \mathcal{R}_{\zeta_{f,s}} \Big( \summ_k
\mathsf{E}_k \hskip1pt T \Delta_{k+s} + \summ_k (id -
\mathsf{E}_k) \hskip1pt T_{4 \cdot 2^{-k}} \Delta_{k+s} \Big)
(f).$$ Assume that $T^*1 = 0$. Our shifted quasi-orthogonal
decomposition in Appendix A asserts that the operator inside the
brackets has norm in $\mathcal{B}(L_2,L_2(\mathcal{H}))$
controlled by $\mathrm{c}_{n,\gamma} s \hskip1pt 2^{- \gamma
s/2}$. In particular, the same happens when we tensor with the
identity on $L_2(\M)$, which is the case. When $T^*1 \neq 0$, we
may follow verbatim the paraproduct argument given in Appendix A
by noting that $\zeta_{f,s} q_k^\perp = q_k^\perp \zeta_{f,s} =
0$. \fin

\begin{remark}
\emph{The projections $$(\1_\Mn - q_k,\zeta_{f,s},\zeta_k)$$
represent the sets $(\Omega_k, \R^n \setminus \Sigma_{f,s}, \R^n
\setminus 9 \hskip1pt \Omega_k)$ in the commutative formulation.}
\end{remark}

\subsection{Weak type $(1,1)$ boundedness}

We prove Theorem B1 in this section. As usual, we may take $f \in
\Mn_{c,+}$. To provide the decomposition $T \hskip-1pt f = A
\hskip-1pt f + B \hskip-1pt f$, we use the projections
$\zeta(\lambda)$ in Lemma \ref{keylem}. Define $$\pi_k =
\bigwedge_{s \ge k} \zeta(2^s) - \bigwedge_{s \ge k-1} \zeta(2^s)
\quad \mbox{for} \quad k \in \Z.$$ If we set $\lambda_f =
\|f\|_\infty$ (which is finite since $f \in \Mn_{c,+}$), we have
$\|f_k\|_\infty \le \lambda_f$ for all integers $k$. In
particular, given any $\lambda > \lambda_f$, it is easily checked
that $q_k(\lambda) = \mathbf{1}_\Mn$ for all $k \in \Z$ and we
deduce that $\zeta(\lambda) = \mathbf{1}_\Mn$ for all $\lambda >
\lambda_f$. This gives rise to $$\sum_{k \in \Z} \pi_k =
\lim_{k_1,k_2 \to \infty} \Big[ \bigwedge_{s \ge k_1} \zeta(2^s) -
\bigwedge_{s \ge -k_2} \zeta(2^s) \Big] = \mathbf{1}_\Mn -
\bigwedge_{k \in \Z} \zeta(2^k) = \mathbf{1}_\Mn - \psi.$$ Our
decomposition for $T \hskip-1pt f$ is the following
\begin{eqnarray*}
T \hskip-1pt f & = & \psi \, T \hskip-1pt f \, \psi \\ [10pt] & +
& (\mathbf{1}_\Mn - \psi) \, T \hskip-1pt f \, \psi + \sum_{j \le
i} \pi_i \, T \hskip-1pt f \, \pi_j \\ & + & \psi \, T \hskip-1pt
f \, (\mathbf{1}_\Mn - \psi) + \sum_{i < j} \pi_i \, T \hskip-1pt
f \, \pi_j \hskip5pt = \hskip5pt \psi \, T \hskip-1pt f \, \psi +
A \hskip-1pt f + B \hskip-1pt f.
\end{eqnarray*}
We claim that $\|Af\|_{L_{1,\infty}(\Mn; \mathcal{H}_r)} + \|B
\hskip-1pt f\|_{L_{1,\infty}(\Mn; \mathcal{H}_c)} \lesssim \,
\|f\|_1$. As we shall see at the end of the proof, the term $\psi
\, T \hskip-1pt f \, \psi$ is even easier to handle. Assume for
simplicity that $\mathcal{H}$ is separable and fix an orthonormal
basis $(u_m)_{m \ge 1}$ of $\mathcal{H}$. If $k_m(x,y) = \langle
u_m, k(x,y) \rangle$, we denote by $T_m \hskip-1pt f$ the
Calder\'on-Zygmund operator associated to the kernel $k_m$, while
$A_m \hskip-1pt f$ and $B_m \hskip-1pt f$ stand for the
corresponding parts of $T_m \hskip-1pt f$. Then, we have $$\big\|
Af \big\|_{L_{1,\infty}(\Mn; \mathcal{H}_r)} = \sup_{\lambda > 0}
\ \lambda \, \varphi \Big\{ \Big( \sum_{m=1}^\infty (A_m
\hskip-1pt f)(A_m \hskip-1pt f)^*\Big)^\frac12 > \lambda \Big\}.$$
As in the martingale case, we may assume that $\lambda = 2^\ell$
for some integer $\ell$. Note that we do not impose $\ell \ge 0$
since $\Mn$ is no longer finite and we also need to consider the
case $0 < \lambda < 1$. Let us define $$w_\ell = \bigwedge_{s \ge
\ell} \zeta(2^s).$$ By the quasi-triangle inequality, we majorize
again by
$$\lambda \, \varphi \Big\{ w_\ell \Big( \sum_{m=1}^\infty (A_m
\hskip-1pt f) (A_m \hskip-1pt f)^* \Big) w_\ell > \lambda^2 \Big\}
+ \lambda \, \varphi \big( \mathbf{1}_\Mn - w_\ell \big) =
\mathsf{A}_1 + \mathsf{A}_2.$$ Then we apply Lemma \ref{keylem} to
estimate $\mathsf{A}_2$ and obtain $$\mathsf{A}_2 \le 2^\ell
\sum_{s \ge \ell} \varphi(\mathbf{1}_\Mn - \zeta(2^s)) \le 9^n \,
2^\ell \sum_{s \ge \ell} 2^{-s} \|f\|_1 \le 2 \cdot 9^n \,
\|f\|_1.$$ We will use that $\psi w_\ell = w_\ell \psi = \psi$ and
$$\pi_k w_\ell = w_\ell \pi_k =
\begin{cases} \pi_k & \mbox{if} \ k \le \ell, \\ 0 &
\mbox{otherwise}. \end{cases}$$ Therefore, defining $\rho_k = \psi
+ \sum_{j \le k} \pi_j$, we have
$$w_\ell A_m \hskip-1pt f = \sum_{i \le
\ell} \pi_i \, (w_\ell \, T_m \hskip-1pt f \, w_\ell) \, \rho_i =
A_{m\ell} \hskip-1pt f.$$ Note that $w_{\ell} T_m \hskip-1pt f
\neq w_{\ell} T_m \hskip-1pt f w_{\ell}$ and that is where we need
to break $T_mf$ into row/column terms. On the other hand, by the
Calder\'on-Zygmund decomposition of $(f,\lambda)$, we conclude
\begin{eqnarray*}
\mathsf{A}_1 & \lesssim & \lambda \, \varphi \Big\{
\sum_{m=1}^\infty (A_{m \ell} g_d) (A_{m \ell} g_d)^* > \lambda^2
\Big\} \\ & + & \lambda \, \varphi \Big\{ \sum_{m=1}^\infty (A_{m
\ell} b_d) (A_{m \ell} b_d)^* > \lambda^2 \Big\} \\ & + & \lambda
\, \varphi \Big\{ \sum_{m=1}^\infty (A_{m \ell} g_{\mathit{off}})
(A_{m \ell} g_{\mathit{off}})^* > \lambda^2 \Big\} \\ & + &
\lambda \, \varphi \Big\{ \sum_{m=1}^\infty (A_{m \ell}
b_{\mathit{off}}) (A_{m \ell} b_{\mathit{off}})^* > \lambda^2
\Big\} \ = \ \mathsf{A}_{g,d} + \mathsf{A}_{b,d} +
\mathsf{A}_{g,\mathit{off}} + \mathsf{A}_{b,\mathit{off}}.
\end{eqnarray*}

\noindent \textbf{The term $\mathsf{A}_{g,d}$.} Chebychev's
inequality gives for $\mathsf{A}_{g,d}$
\begin{eqnarray*}
\mathsf{A}_{g,d} & \le & \frac{1}{\lambda} \, \sum_{m=1}^\infty
\sum_{i \le \ell} \varphi \big( \pi_i \, (w_\ell \, T_m g_d \,
w_\ell) \, \rho_i (w_\ell \, T_m g_d \, w_\ell)^* \, \pi_i \big)
\\ & \le & \frac{1}{\lambda} \, \sum_{m=1}^\infty \varphi \big(
(w_\ell \, T_m g_d \, w_\ell) \, (w_\ell \, T_m g_d \, w_\ell)^*
\big) \\ & \le & \frac{1}{\lambda} \, \sum_{m=1}^\infty \varphi
\big( (T_m g_d) \, (T_m g_d)^* \big) \, = \, \frac{1}{\lambda} \,
\big\| T g_d \big\|_{L_2(\Mn; \mathcal{H}_{oh})}^2 \, \le \,
\frac{1}{\lambda} \, \|g_d\|_2^2 \, \le \, 2^n \|f\|_1.
\end{eqnarray*}
Last inequality is part of the statement of Calder\'on-Zygmund
decomposition above.

\noindent \textbf{The term $\mathsf{A}_{g,\mathit{off}}$.} Arguing
as for $g_d$ we find
\begin{eqnarray*}
\mathsf{A}_{g,\mathit{off}} & \le & \frac{1}{\lambda} \,
\sum_{m=1}^\infty \varphi \big( (w_\ell \, T_m g_{\mathit{off}} \,
w_\ell) \, (w_\ell \, T_m g_{\mathit{off}} \, w_\ell)^* \big) \\ &
= & \frac{1}{\lambda} \int_{\R^n} \big\| \big[ w_\ell \, T
g_{\mathit{off}} \, w_\ell \big] (x)
\big\|_{L_2(\M;\mathcal{H}_{oh})}^2 \, dx.
\end{eqnarray*}
On the other hand, we have $$\mathrm{supp}^* {dg_{(s)}}_{k+s} \le
\mathbf{1}_\Mn - q_k(2^\ell)$$ from Remark \ref{goff} and we claim
$w_\ell \le \zeta_{g_{(s)},s}$. Indeed, clearly $w_\ell \le
\zeta(2^\ell)$ and
\begin{eqnarray*}
\bigvee_{Q \in \mathcal{Q}_k} (\mathbf{1}_\M - \xi_Q) 1_{9Q} & = &
\mathrm{supp} \ \sum_{Q \in \mathcal{Q}_k} (\mathbf{1}_\M - \xi_Q)
1_{9Q} \\ & \le & \mathrm{supp} \sum_{s=m_\lambda+1}^k \sum_{Q \in
\mathcal{Q}_s} (\xi_{\widehat{Q}} - \xi_Q) 1_{9Q} \ = \
\mathrm{supp} \, \psi_k(2^\ell).
\end{eqnarray*}
Therefore, we conclude $$w_\ell \le \zeta(2^\ell) = \bigwedge_{k
\in \Z} \big( \mathbf{1}_\Mn - \mathrm{supp} \, \psi_k(2^\ell)
\big) \le \bigwedge_{k \in \Z} \big( \mathbf{1}_\Mn - \bigvee_{Q
\in \mathcal{Q}_k} (\mathbf{1}_\M - \xi_Q) 1_{9Q} \big) =
\zeta_{g_{(s)},s}.$$ Using this inequality and
pseudo-localization (Theorem \ref{pseudolocal}), we obtain
\begin{eqnarray*}
\mathsf{A}_{g,\mathit{off}} & \le & \frac{1}{\lambda} \int_{\R^n}
\big\| \big[ \zeta_{g_{(s)},s} \, T g_{\mathit{off}} \,
\zeta_{g_{(s)},s} \big] (x) \big\|_{L_2(\M;\mathcal{H}_{oh})}^2 \,
dx \\ & \le & \frac{1}{\lambda} \, \Big( \sum_{s=1}^\infty \Big[
\int_{\R^n} \big\| \big[ \zeta_{g_{(s)},s} \, T g_{(s)} \,
\zeta_{g_{(s)},s} \big] (x) \big\|_{L_2(\M;\mathcal{H}_{oh})}^2 \,
dx \Big]^\frac12 \Big)^2 \\ & \le &
\frac{\mathrm{c}_{n,\gamma}}{\lambda} \, \Big( \sum_{s=1}^\infty s
\, e^{-\gamma s/2} \, \Big[ \sum_{k=m_\lambda+1}^\infty
\|g_{k,s}\|_2^2 \Big]^\frac12 \Big)^2 \ \lesssim \
\mathrm{c}_{n,\gamma} \|f\|_1,
\end{eqnarray*}
where the last estimate follows from the inequality given in
Remark \ref{goff}.

\noindent \textbf{The term $\mathsf{A}_{b,d}$.} We have
$$A_{m\ell} b_d = \sum_{j \le i \le \ell} \pi_i \, (w_\ell \, T_m
\hskip-1pt b_d \, w_\ell) \, \pi_j + \Big[ \sum_{i \le \ell} \pi_i
\Big] \, (w_\ell \, T_m \hskip-1pt b_d \, w_\ell) \, \psi.$$ Since
$\sum_{j \le i \le \ell} \pi_i \, \cdot \, \pi_j$ and $[\sum_{i
\le \ell} \pi_i] \, \cdot \, \psi$ are of weak type $(1,1)$
$$\mathsf{A}_{b,d} \le \Big\| \sum_{m=1}^\infty A_{m \ell} b_d
\otimes e_{1m} \Big\|_{1,\infty} \lesssim \Big\| \sum_{m=1}^\infty
(w_\ell \, T_m \hskip-1pt b_d \, w_\ell) \otimes e_{1m}
\Big\|_1.$$ Letting $b_{d,k} = p_k (f-f_k) p_k$, we find that
$\mathsf{A}_{b,d}$ is dominated by
$$\sum_{k=m_\lambda+1}^\infty
\sum_{Q\in \mathcal{Q}_k} \tau \otimes \int_{\R^n} \Big(
\sum_{m=1}^\infty \Big| \big[ \int_Q k_m(x,y) (w_\ell(x) \,
b_{d,k}(y) \, w_\ell(x)) \, dy \big]^* \Big|^2 \Big)^\frac12 \,
dx.$$ Given $y \in Q \in \mathcal{Q}_k$, we have $p_k(y) =
\xi_{\widehat{Q}} - \xi_Q$ and $$w_\ell(x) \, b_{d,k}(y) \,
w_\ell(x) = 0 \quad \mbox{for} \quad (x,y) \in 9Q \times Q$$ by
Lemma \ref{keylem}. This means that
\begin{eqnarray*}
\mathsf{A}_{b,d} & \lesssim & \sum_{k=m_\lambda+1}^\infty
\sum_{Q\in \mathcal{Q}_k} \Big\| w_\ell \Big( \sum_{m=1}^\infty
1_{\R^n \setminus 9Q} \, T_m \hskip-1pt (b_{d,k} 1_Q) \otimes
e_{1m} \Big) w_\ell \Big\|_1 \\ & \le &
\sum_{k=m_\lambda+1}^\infty \sum_{Q \in \mathcal{Q}_k} \tau
\otimes \int_{ \mathbb{R}^n \setminus 9Q} \Big( \sum_{m=1}^\infty
\Big| [ \int_Q k_m(x,y) b_{d,k}(y) \, dy]^* \Big|^2
\Big)^{\frac12} \, dx.
\end{eqnarray*}
By the mean zero of $b_{d,k}$ on $Q$, we rewrite our term as
follows $$\sum_{k=m_\lambda+1}^\infty \sum_{Q \in \mathcal{Q}_k}
\tau \otimes \int_{ \mathbb{R}^n \setminus 9Q} \Big(
\sum_{m=1}^\infty \Big| [ \int_Q \big( k_m(x,y)-k_m(x,c_Q) \big)
\, b_{d,k}(y) \, dy]^* \Big|^2 \Big)^{\frac12} \, dx,$$ where
$c_Q$ denotes the center of $Q$. We see it as an $L_1(\mathcal{A};
\mathcal{H}_r)$ norm of $$\sum_{m=1}^\infty \Big( \int_Q \big(
k_m(x,y) - k_m(x,c_Q) \big) \, b_{d,k}(y) \, dy \Big) \otimes
u_m.$$ Putting the $L_1(\mathcal{A}, \mathcal{H}_r)$ norm into the
integral $\int_Q$ gives us a larger term
\begin{eqnarray*}
\lefteqn{\sum_{k=m_\lambda + 1}^\infty \sum_{Q \in \mathcal{Q}_k}
\tau \otimes \int_{\R^n \setminus 9Q} \int_Q \Big(
\sum_{m=1}^\infty \Big| b_{d,k}^{*}(y) \, \big( k_m^{*}(x,y) -
k_m^{*}(x,c_Q) \big) \Big|^2 \Big)^{\frac12} \, dy \, dx} \\
& = & \sum_{k=m_\lambda+1}^\infty \sum_{Q \in \mathcal{Q}_k} \tau
\otimes \int_{\R^n \setminus 9Q} \int_Q \Big( \sum_{m=1}^\infty
\big| k_m^{*}(x,y) - k_m^{*}(x,c_Q) \big|^2 \Big)^{\frac12} \,
|b_{d,k}^{*}(y)| \, dy \, dx \\ & \lesssim &
\sum_{k=m_\lambda+1}^\infty \sum_{Q \in \mathcal{Q}_k} \Big[
\int_{\R^n \setminus 9Q}
\frac{|\ell(Q)|^\gamma}{|x-c_Q|^{n+\gamma}} \, dx \Big] \, \Big[
\tau \otimes \int_Q |b_{d,k}^{*}(y)| \, dy \Big] \\
& \lesssim & \sum_{k=m_\lambda+1}^\infty \tau \otimes \int_{\R^n}
|b_{d,k}^{*}(y)| \, dy \ \le \ 2 \, \|f\|_1.
\end{eqnarray*}
The last inequality uses the diagonal estimate in
Calder\'on-Zygmund decomposition.

\noindent \textbf{The term $\mathsf{A}_{b,\mathit{off}}$.}
Although it is more technical, the estimate for
$\mathsf{A}_{b,\mathit{off}}$ is very similar in nature to that
for $\mathsf{A}_{b,d}$. Indeed, the only significant difference
relies on the fact that we have to show
$$\sum_{k=m_\lambda+1}^\infty \Big\| \sum_{m=1}^\infty (w_\ell \, T_m
b_{k,s} \, w_\ell) \otimes e_{1m} \Big\|_1 \lesssim 2^{-\gamma s}
\, \|f\|_1$$ for each $s \ge 1$. Namely, here we write $b_{k,s}$
for the sum of the $k$-th terms in the upper and lower $s$-th
diagonals of $b_{\mathit{off}}$, and the estimate above resembles
the same procedure that we used for $g_{\mathit{off}}$ in
comparison with $g_d$. That is, there exists a geometric
\emph{almost diagonal} phenomenon. It is nevertheless
straightforward to check that the arguments above for
$\mathsf{A}_{b,d}$ also apply for $\mathsf{A}_{b,\mathit{off}}$
using the ideas in \cite{P2} for the off-diagonal term of the bad
part in Calder\'on-Zygmund decomposition.

\noindent \textbf{Conclusion.} By symmetry, the same arguments
apply to estimate the norm of the column term $B \hskip-1pt f$. It
therefore remains to consider the term $\psi \, T \hskip-1pt f \,
\psi$. However, note that $\psi \le w_\ell$ and we know from the
arguments above how to estimate the term $w_\ell \, T \hskip-1pt f
\, w_\ell$. The proof is completed. \fin

\subsection{$\mathrm{BMO}$ estimate, interpolation and duality}

We now prove Theorem B2. Let us briefly recall the definition of
the noncommutative analogue of function-BMO spaces. According to
\cite{Me}, we define the norms (modulo constants)
\begin{eqnarray*}
\|f\|_{\mathrm{BMO}_\Mn^r} & = & \sup_{\mathcal{Q} \ \mathrm{cube}
\, \subset \R^n} \Big\| \Big( \frac{1}{|Q|} \int_Q
(f(x)-f_Q)(f(x)-f_Q)^* \, dx \Big)^{\frac12} \Big\|_\M, \\
\|f\|_{\mathrm{BMO}_\Mn^c} & = & \sup_{\mathcal{Q} \ \mathrm{cube}
\, \subset \R^n} \Big\| \Big( \frac{1}{|Q|} \int_Q
(f(x)-f_Q)^*(f(x)-f_Q) \, dx \Big)^{\frac12} \Big\|_\M.
\end{eqnarray*}
We also set $$\|f\| _{\mathrm{BMO}(\Mn)} = \max \Big\{
\|f\|_{\mathrm{BMO}_\Mn^r}, \|f\|_{\mathrm{BMO}_\Mn^c} \Big\}.$$
The spaces $\mathrm{BMO}(\Mn;\mathcal{H}_r)$ and
$\mathrm{BMO}(\Mn;\mathcal{H}_c)$ were defined in the
Introduction. As in the previous paragraph, we fix an orthonormal
basis $(u_m)_{m \ge 1}$ of $\mathcal{H}$ and define $k_m(x,y)$ and
$T_m \hskip-1pt f$ accordingly. Letting $\mathcal{R} = \Mn
\bar\otimes \mathcal{B}(\ell_2)$, we have
\begin{eqnarray*}
\|T \hskip-1pt f\|_{\mathrm{BMO}(\Mn; \mathcal{H}_r)} & = & \Big\|
\sum_{m=1}^\infty T_m \hskip-1pt f \otimes e_{1m}
\Big\|_{\mathrm{BMO}(\mathcal{R})}, \\ \|T \hskip-1pt
f\|_{\mathrm{BMO}(\Mn; \mathcal{H}_c)} & = & \Big\|
\sum_{m=1}^\infty T_m \hskip-1pt f \otimes e_{m1}
\Big\|_{\mathrm{BMO}(\mathcal{R})}.
\end{eqnarray*}
By symmetry, we just prove $\|T \hskip-1pt
f\|_{\mathrm{BMO}(\Mn;\mathcal{H}_c)} \lesssim \|f\|_\infty$. As
usual, in the definition of the $\mathrm{BMO}$ norm of a function
$f$, we may replace the averages $f_Q$ by any other operator
$\alpha_Q$ depending on $Q$. Fix a cube $Q$ in $\R^n$ and let
$$\alpha _Q = \sum_{m=1}^\infty \alpha_{Q,m} \otimes e_{m1}
= \sum_{m=1}^\infty \frac 1{|Q|} \int_Q \Big( \int_{\R^n \setminus
2Q} k_m(z,y)f(y) \, dy \Big) \, dz \otimes e_{m1}.$$ We have
\begin{eqnarray*}
\lefteqn{\hskip-10pt \big[ T_m \hskip-1pt f - \alpha_{Q,m}
\big](x)} \\ & = & \frac{1}{|Q|} \int_Q \int_{\R^n \setminus 2Q}
\big( k_m(x,y) - k_m(z,y) \big) f(y) \, dy \, dz + \int_{2Q}
k_m(x,y) f(y) \, dy
\\ & = & B_{m1} \hskip-1pt f(x) + B_{m2} \hskip-1pt f(x).
\end{eqnarray*}
For $B_1 \hskip-1pt f = \sum_m B_{m1} \hskip-1pt f \otimes
e_{m1}$, we have
\begin{eqnarray*}
\lefteqn{\hskip-20pt \sup_{x \in Q} \|B_1 \hskip-1pt f(x)\|_{\M
\bar\otimes \mathcal{B}(\ell_2)}} \\ & \le & \sup_{x,z \in Q}
\Big\| \sum_{m=1}^\infty \int_{\R^n \setminus 2Q} \big( k_m(x,y) -
k_m(z,y) \big) f(y) \, dy \otimes e_{m1} \Big\|_{\M \bar\otimes
\mathcal{B}(\ell_2)} \\ & = & \sup_{x,z \in Q} \Big\|
\sum_{m=1}^\infty \Big| \int_{\R^n \setminus 2Q} \big( k_m(x,y) -
k_m(z,y) \big) f(y) \, dy \Big|^2 \Big\|_\M^{\frac12} \\ & \le &
\sup_{x,z \in Q} \Big( \sum_{m=1}^\infty \big[ \int_{\R^n
\setminus 2Q} \big| k_m(x,y) - k_m(z,y) \big| \, dy \big]^2
\Big)^\frac12 \, \|f\|_\infty \\ & \le & \sup_{x,z \in Q} \Big[
\int_{\R^n \setminus 2Q} \Big( \sum_{m=1}^\infty \big| k_m(x,y) -
k_m(z,y) \big|^2 \Big)^\frac12 dy \Big] \, \|f\|_\infty \\ & \le &
\sup_{x,z \in Q} \Big[ \int_{\R^n \setminus 2Q}
\frac{\ell(Q)^\gamma}{|x-y|^{n+\gamma}} \, dy \Big] \,
\|f\|_\infty \ \lesssim \ \|f\|_\infty.
\end{eqnarray*}
For $B_2 \hskip-1pt f = \sum_m B_{m2} \hskip-1pt f \otimes
e_{m1}$, we have
\begin{eqnarray*}
\lefteqn{\hskip-20pt \frac{1}{|Q|} \Big\| \int_Q (B_2 \hskip-1pt
f(x))^* (B_2 \hskip-1pt f(x)) \, dx \Big\|_{\M \bar\otimes
\mathcal{B}(\ell_2)}} \\ & = & \frac{1}{|Q|} \Big\| \int_Q \Big|
\sum_{m=1}^\infty \int_{2Q} k_m(x,y) f(y) \, dy \otimes e_{m1}
\Big|^2 \, dx \Big\|_{\M \bar\otimes \mathcal{B}(\ell_2)} \\
& = & \frac{1}{|Q|} \Big\| \int_Q \sum_{m=1}^\infty \Big|
\int_{2Q} k_m(x,y) f(y) \, dy \Big|^2 \, dx \Big\|_\M \\ & = &
\frac{1}{|Q|} \sup_{\|a\|_{L_2(\M)} \le 1} \tau \otimes \int_Q
\hskip4pt \sum_{m=1}^\infty \Big| \int_{2Q} k_m(x,y) f(y) a \, dy
\Big|^2 dx \\ & \le & \frac{1}{|Q|} \sup_{\|a\|_{L_2(\M)} \le 1}
\tau \otimes \int_{\R^n} \sum_{m=1}^\infty \Big|
\int_{\R^n} k_m(x,y) f(y) a 1_{2Q}(y) \, dy \Big|^2 dx \\
& = & \frac{1}{|Q|} \sup_{\|a\|_{L_2(\M)}
\le 1} \big\| T(fa1_{2Q}) \big\|_{L_2(\Mn; \mathcal{H}_{oh})}^2 \\
& \lesssim & \frac{1}{|Q|} \sup_{\|a\|_{L_2(\M)} \le 1} \|fa
1_{2Q}\|_{L_2(\Mn)}^2 \ \lesssim \ \|f\|_\infty^2.
\end{eqnarray*}
On the other hand, if $\beta_mf(x) = \int_{2Q} k_m(x,y) f(y) \,
dy$, we also have
\begin{eqnarray*}
\lefteqn{\hskip-20pt \frac{1}{|Q|} \Big\| \int_Q (B_2 \hskip-1pt
f(x))(B_2 \hskip-1pt f(x))^* \, dx \Big\|_{\M \bar\otimes
\mathcal{B}(\ell_2)}} \\ & = & \frac{1}{|Q|} \Big\|
\sum_{m_1,m_2=1}^\infty \Big[ \int_Q \beta_{m_1}f(x)
\beta_{m_2}^*f(x) \, dx \Big] \otimes e_{m_1,m_2} \Big\|_{\M
\bar\otimes \mathcal{B}(\ell_2)} \\ & = & \frac{1}{|Q|}
\|\Lambda\|_{\M \bar\otimes \mathcal{B}(\ell_2)}.
\end{eqnarray*}
Since $\Lambda$ is a positive operator acting on $L_2(\M; \ell_2)$
\begin{eqnarray*}
\lefteqn{\hskip-10pt \frac{1}{|Q|} \Big\| \int_Q (B_2 \hskip-1pt
f(x))(B_2 \hskip-1pt f(x))^* \, dx \Big\|_{\M \bar\otimes
\mathcal{B}(\ell_2)}} \\ & = & \frac{1}{|Q|} \sup_{\|a\|_{L_2(\M;
\ell_2)} \le 1} \big\langle \Lambda a, a \big\rangle \\ & = &
\frac{1}{|Q|} \sup_{\|a\|_{L_2(\M; \ell_2)} \le 1} \hskip4.7pt
\tau \Big( \big[ \sum_{m_1=1}^\infty a_{m_1}^* \otimes e_{1m_1}
\big] \, \Lambda \,
\big[ \sum_{m_2=1}^\infty a_{m_2} \otimes e_{m_21} \big] \Big) \\
& = & \frac{1}{|Q|} \sup_{\|a\|_{L_2(\M; \ell_2)} \le 1}
\hskip4.7pt \tau \otimes \int_Q \Big| \sum_{m=1}^\infty
\beta_m^*f(x) a_m \Big|^2 \, dx \\ & = & \frac{1}{|Q|}
\sup_{\|a\|_{L_2(\M; \ell_2)} \le 1} \hskip4.7pt \tau \otimes
\int_Q \Big| \sum_{m=1}^\infty \int_{\R^n} \overline{k_m(x,y)}
f(y)^* a_m 1_{2Q}(y) \, dy \Big|^2 \, dx \\ & \le & \frac{1}{|Q|}
\sup_{\begin{subarray}{c} \|a\|_{L_2(\M; \ell_2)} \le 1 \\
\|g\|_{L_2(\Mn)} \le 1
\end{subarray}} \Big[ \tau \otimes \int_{\R^n} \sum_{m=1}^\infty
\int_{\R^n} \overline{k_m(x,y)} f(y)^* a_m 1_{2Q}(y) \, dy \, g(x)
\, dx \Big]^2 \\ & = & \frac{1}{|Q|} \sup_{\begin{subarray}{c}
\|a\|_{L_2(\M; \ell_2)} \le 1 \\ \|g\|_{L_2(\Mn)} \le 1
\end{subarray}} \Big[ \tau \otimes \int_{\R^n}
\sum_{m=1}^\infty \Big( \int_{\R^n} \overline{k_m(x,y)} g(x) \, dx
\Big) f(y)^* a_m 1_{2Q}(y) \, dy \Big]^2 \\ & \le & \frac{1}{|Q|}
\sup_{\begin{subarray}{c} \|a\|_{L_2(\M; \ell_2)} \le 1 \\
\|g\|_{L_2(\Mn)} \le 1 \end{subarray}} \|T \hskip-1pt
g^*\|_{L_2(\Mn; \mathcal{H}_{oh})}^2 \, \Big( \tau \otimes
\int_{\R^n} \sum_{m=1}^\infty |f(y)^* a_m|^2 1_{2Q}(y) \, dy \Big)
\\ & \le & \frac{1}{|Q|}
\sup_{\begin{subarray}{c} \|a\|_{L_2(\M; \ell_2)} \le 1 \\
\|g\|_{L_2(\Mn)} \le 1 \end{subarray}} \|T \hskip-1pt
g^*\|_{L_2(\Mn; \mathcal{H}_{oh})}^2 \, \tau \Big(
\sum_{m=1}^\infty |a_m|^2 \Big) |2Q| \, \|f\|_\infty^2 \ \lesssim
\ \|f\|_\infty^2.
\end{eqnarray*}
The estimates so far and their row analogues give rise to
$$\max \Big\{ \big\| T \hskip-1pt f \big\|_{\mathrm{BMO}(\Mn;
\mathcal{H}_r)}, \big\| T \hskip-1pt f \big\|_{\mathrm{BMO}(\Mn;
\mathcal{H}_c)} \Big\} \lesssim \|f\|_\infty.$$ This completes the
BMO estimate. By interpolation as in Section \ref{S1}, we get
$$\big\| T \hskip-1pt f \big\|_{L_p(\Mn; \mathcal{H}_{rc})}
\le c_p \, \|f\|_p$$ for $1 < p < \infty$. Finally, if $T$ is an
isometry $L_2(\Mn) \to L_2(\Mn;\mathcal{H}_{oh})$, polarization
gives \\ [8pt] \null \hfill $\displaystyle \|f\|_{L_p(\Mn)} =
\sup_{\|g\|_{L_{p'}(\Mn)} \le 1} \langle Tf,Tg \rangle \lesssim
\big\| T \hskip-1pt f \big\|_{L_p(\Mn; \mathcal{H}_{rc})}.$ \hfill
$\square$

\subsection{Examples and further comments}

The first examples that come to mind are the scalar-valued
Calder\'on-Zygmund operators studied in \cite{P2}, given by
$\mathcal{H}$ one dimensional. On the other hand, as in \cite{P2}
or in Section \ref{S1} above, we might also consider commuting
operator-valued coefficients in Theorems B1 and B2. We omit the
formal statement of such result. Let us study some more
examples/applications of our results.

\noindent A. \emph{Lusin square functions and $g$-functions.} A
row/column form of the Lusin area function and the
Littlewood-Paley $g$-function for the Poisson kernel was given in
\cite{Me}. It is standard that both square functions are
associated to Calder\'on-Zygmund kernels satisfying our
size/smoothness conditions. Moreover, the $L_2$ boundedness of
these operators is straightforward. Therefore, Theorems B1 and B2
apply and we obtain the weak $L_1$, strong $L_p$ and
$L_\infty-\mathrm{BMO}$ boundedness of these operators in the
operator valued setting. The strong $L_p$ inequalities for Lusin
square functions and $g$-functions were one of the main results in
\cite{Me}, while the weak $L_1$ and $L_\infty-\mathrm{BMO}$
estimates are new. In a similar way, we may consider any other
symmetric diffusion semigroup as far as it satisfies our kernel
assumptions. Note that, when dealing with strong $L_p$ estimates,
more general families of semigroups were considered in \cite{JLX}.
It seems however that the approach in \cite{JLX} does not permit
to work with general Calder\'on-Zygmund operators (not coming from
semigroups) or providing weak $L_1$ and $L_\infty-\mathrm{BMO}$
estimates.

\noindent B. \emph{Operator-valued Littlewood-Paley theorem (with
better constants).} Consider $$\widehat{M_k \hskip-1pt f} =
1_{\Delta_k} \widehat{f} \quad \mbox{with} \quad \Delta_k =
(-2^{k+1},-2^k] \cup [2^k, 2^{k+1}).$$ An immediate application of
Theorem B2 is a generalization to the operator valued setting of
the Littlewood-Paley theorem, which assets in the commutative case
that for any $1 < p < \infty$ we have
$$\|f\|_p \sim_{c_p} \Big\| \Big( \summ_k |M_k \hskip-1pt f|^2
\Big)^\frac12 \Big\|_p.$$ Indeed, Theorem B2 easily gives the
$L_\infty$-$\mathrm{BMO}$ estimate by the boundedness of Riesz
transforms on BMO and the classical ``shift-truncation" argument
of smoothing multipliers. An ${\mathcal H}_1$-weak (1,1) estimate
can be also obtained by considering atomic decomposition \cite{Ri}
and a dual version of Theorem B2. However, using the fact that
$L_p(\M)$ is a UMD Banach space, this also follows from Bourgain's
vector-valued Littlewood-Paley inequality for UMD spaces
\cite{Bo3}
$$\|f\|_{L_p(\mathbb{X})} \sim_{C_p} \Big( \int_\Omega \Big\|
\summ_k M_k \hskip-1pt f \, r_k(w) \Big\|_{L_p(\mathbb{X})}^p  \,
dw \Big)^\frac1p$$ combined with the noncommutative Khintchine
inequality. In conclusion, Theorem B2 provides a new proof of
Bourgain's result for $L_p(\M)$-valued functions. The advantage is
that we will get the optimal constants $c_p\sim\frac {1}{p-1} , p$
as $p \to 1,\infty. ^{\dag}$ \footnote {\dag\  Musat's constants
were improved to $\sim p$ for $p>2$ after Randrianantoanina's work
\cite{R3}.}

\noindent C. \emph{Beyond the Lebesgue measure on $\R^n$.} Let us
point some possible generalizations of our results for future
research. First, following Han and Sawyer \cite{HS}, we may
consider operator-valued Littlewood-Paley inequalities on
homogeneous spaces. We are confident these results should hold.
Second, in a less obvious way, we might work in the nondoubling
setting. Recall that Littlewood-Paley theory for nondoubling
measures was a corner stone through the $T1$ theorem by Nazarov,
Treil and Volberg \cite{NTV2} and Tolsa \cite{To}. Tolsa's
technique seems reasonable as far as we know how to prove the weak
type (1,1) inequality. The Calder\'on-Zygmund decomposition in
\cite{To2} looks like the most difficult step and an interesting
problem.

\noindent D. \emph{An application to the fully noncommutative
setting.} In a forthcoming paper \cite{JMP}, we are going to apply
our results in this paper to Fourier multipliers on group von
Neumann algebras $VN(G)$. The idea is to embed $VN(G)$ into
$L_\infty({\Bbb R}^n) \bar\otimes \mathcal{B}(\ell_2(G))$ and to
reduce the boundedness of Fourier multipliers on $VN(G)$ to the
boundedness of singular integrals studied here.

\section*{Appendix A. Hilbert space valued
pseudo-localization}

\renewcommand{\theequation}{A\arabic{equation}}

\setcounter{equation}{0}

Let us sketch the modifications from the argument in \cite{P2}
needed to extend the pseudo-localization result there to the
context of Hilbert space valued kernels. We adopt the terminology
from Section \ref{S2}. Besides, we shall write just $L_p$ to refer
to the commutative $L_p$ space on $\R^n$ equipped with the
Lebesgue measure $dx$ and $L_p(\mathcal{H})$ for its
$\mathcal{H}$-valued extension.

\begin{localtheo}
Given a Hilbert space $\mathcal{H}$, let us consider a kernel $k:
\R^{2n} \setminus \Delta \to \mathcal{H}$ satisfying the
size/smoothness conditions imposed in the Introduction, which
formally defines $$T \! f(x) = \int_{\R^n} k(x,y) f(y) \, dy,$$ a
Calder\'on-Zygmund operator. Assume further that $T: L_2 \to
L_2(\mathcal{H})$ is of norm $1$. Let us fix a positive integer
$s$. Given a function $f$ in $L_2$ and any integer $k$, we define
$\Omega_k$ to be the smallest $\mathcal{R}_k$-set containing the
support of $df_{k+s}$. If we further consider the set
$$\Sigma_{f,s} = \bigcup_{k \in \Z} 9 \hskip1pt \Omega_k,$$ then
we have the localization estimate
$$\Big( \int_{\R^n \setminus \hskip1pt
\Sigma_{f \! ,s}} \|T \! f(x)\|_\mathcal{H}^2 \, dx
\Big)^{\frac12} \ \le \ \mathrm{c}_{n,\gamma} \hskip1pt s
\hskip1pt 2^{- \gamma s/2} \Big( \int_{\R^n} |f(x)|^2 dx
\Big)^{\frac12}.$$
\end{localtheo}

\noindent This appendix is devoted to sketch the proof of the
result stated above.

\subsection*{\textnormal{A.1.} Three auxiliary results}
\label{3AR}

As in \cite{P2}, Cotlar and Schur lemmas as well as the
localization estimate from \cite{MC} are the building blocks of
the argument. We shall need some Hilbert space valued forms of
these results. The exact statements are given below. The proofs
are simple generalizations.

\begin{cotlar}
Given Hilbert spaces $\mathcal{K}_1,\mathcal{K}_2$, let us
consider a family $(T_k)_{k \in \Z}$ of bounded operators
$\mathcal{K}_1 \to \mathcal{K}_2$ with finitely many non-zero
$T_k$'s. Assume that there exists a sumable sequence
$(\alpha_k)_{k \in \Z}$ such that
$$\max \Big\{ \big\| T_i^* T_j^{\null}
\big\|_{\mathcal{B}(\mathcal{K}_1)}, \big\| T_i^{\null} T_j^*
\big\|_{\mathcal{B}(\mathcal{K}_2)} \Big\} \, \le \,
\alpha_{i-j}^2$$ for all $i,j \in \Z$. Then we automatically have
$$\Big\| \summ_k T_k \Big\|_{\mathcal{B}(\mathcal{K}_1,
\mathcal{K}_2)} \, \le \, \summ_k \alpha_k.$$
\end{cotlar}

\begin{schur}
Let $T$ be given by $$T \! f(x) = \int_{\R^n} k(x,y) \hskip1pt
f(y) \, dy.$$ Let us define the Schur integrals associated to $k$
$$\mathcal{S}_1(x) = \int_{\R^n} \big\| k(x,y) \big\|_\mathcal{H}
\, dy \quad \mbox{and} \quad \mathcal{S}_2(y) = \int_{\R^n} \big\|
k(x,y) \big\|_\mathcal{H} \, dx.$$ If $\mathcal{S}_1,
\mathcal{S}_2 \in L_\infty$, then $T$ is bounded on $L_2$ and we
have $$\|T\|_{\mathcal{B}(L_2, L_2(\mathcal{H}))} \le \sqrt{
\big\| \mathcal{S}_1 \big\|_\infty \big\| \mathcal{S}_2
\big\|_\infty^{\null} }.$$
\end{schur}

\begin{localest} Assume that $$\big\| k(x,y) \big\|_\mathcal{H}
\lesssim \frac{1}{|x-y|^n} \quad \mbox{for all} \quad x,y \in
\R^n,$$ and let $T$ be a Calder{\'o}n-Zygmund operator associated
to the kernel $k$. Assume further that $T: L_2 \to
L_2(\mathcal{H})$ is of norm $1$. Then, given $x_0 \in \R^n$ and
$r_1, r_2 \in \R_+$ with $r_2 > 2 \hskip1pt r_1$, the estimate
below holds for any pair $f,g$ of bounded scalar-valued functions
respectively supported by $\mathsf{B}_{r_1}(x_0)$ and
$\mathsf{B}_{r_2}(x_0)$
$$\Big\| \int_{\R^n} T \! f(x) g(x) \, dx \Big\|_\mathcal{H} \le
\mathrm{c}_n \hskip1pt r_1^n \hskip1pt \log(r_2/r_1) \hskip1pt
\|f\|_\infty \|g\|_\infty.$$
\end{localest}

\subsection*{\textnormal{A.2.} A quasi-orthogonal decomposition}

According to the conditions imposed on $T$, it is clear that its
adjoint $T^*: L_2(\mathcal{H}^*) \to L_2$ is a norm $1$ operator
with kernel given by $k^*(x,y) = \langle k(y,x), \cdot \rangle$.
Indeed, we have
\begin{eqnarray*}
\big\langle T \hskip-1pt f,g \big\rangle & = & \int_{\R^n}
\Big\langle \int_{\R^n} k(x,y) f(y) \, dy , g(x) \Big\rangle \, dx
\\ & = & \int_{\R^n} \overline{f(y)} \int_{\R^n} \langle k(x,y),
g(x) \rangle \, dx \, dy \ = \ \big\langle f, T^*g \big\rangle,
\end{eqnarray*}
for any $g \in L_2(\mathcal{H}^*)$. In particular, the kernel
$k^*(x,y)$ may be regarded as a linear functional $\mathcal{H}^*
\to \C$ or, equivalently, an element in $\mathcal{H}$. Thus, it
satisfies the same size and smoothness estimates as $k(x,y)$ and
we may and shall view $T^*$ as an operator $L_2 \to
L_2(\mathcal{H})$ which formally maps $g \in L_2$ to
$$T^*g(x) = \int_{\R^n} k^*(x,y) \hskip1pt g(y) \, dy \in
L_2(\mathcal{H}).$$ We shall also use the terminology $\langle T
\hskip-1pt f, g \rangle = \langle f, T^*g \rangle$ to denote the
continuous bilinear form $(f,g) \in L_2 \times L_2 \to
\mathcal{H}$. In fact, we may also define $T^*1$ in a weak sense
as in the scalar-valued case. Moreover, the condition $T^*1=0$
which we shall assume at some points in this section, implies as
usual that the relation below holds for any $f \in H_1$
\begin{equation} \label{T*1=0}
\int_{\R^n}^{\null} T f (x) \, dx = 0.
\end{equation}
Indeed, we formally have $\langle T \hskip-1pt f,1 \rangle =
\langle f, T^*1 \rangle = 0$. Nevertheless, we refer to Hyt\"onen
and Weis \cite{HW2} for a more in depth explanation of all the
identifications we have done so far. Let us go back to our
problem. As usual, let $\mathsf{E}_k$ be the $k$-th dyadic
conditional expectation and fix $\Delta_k$ for the martingale
difference $\mathsf{E}_k - \mathsf{E}_{k-1}$, so that
$\mathsf{E}_k(f) = f_k$ and $\Delta_k(f) = df_k$. Then we consider
the following decomposition $$1_{\R^n \setminus \Sigma_{f,s}} T \!
f = 1_{\R^n \setminus \Sigma_{f,s}} \Big( \sum_{k \in \Z}
\mathsf{E}_k T \Delta_{k+s} + \sum_{k \in \Z} (id-\mathsf{E}_k)
T_{4 \cdot 2^{-k}} \Delta_{k+s} \Big) (f),$$ where
$T_{\varepsilon}$ denotes the truncated singular integral
$$T_\varepsilon f (x) = \int_{|x-y| > \varepsilon} k(x,y)
\hskip1pt f(y) \, dy.$$ We refer to \cite{P2} more details. Our
first step towards the proof is the following.

\begin{shifttheo}
Let $T: L_2 \to L_2(\mathcal{H})$ be a normalized
Calder{\'o}n-Zygmund operator with Lipschitz parameter $\gamma$,
as defined above. Assume further that $T^*1 = 0$, so that
$$\int_{\R^n}^{\null} T \! f(x) \, dx = 0$$ for any $f \in H_1$.
Then, we have
$$\|\Phi_s\|_{\mathcal{B}(L_2,L_2(\mathcal{H}))} = \Big\| \summ_k
\mathsf{E}_k T \Delta_{k+s} \Big\|_{\mathcal{B}(L_2,
L_2(\mathcal{H}))} \le \mathrm{c}_{n,\gamma} \hskip1pt s \hskip1pt
2^{- \gamma s/2}.$$ Moreover, regardless the value of $T^*1$, we
also have
$$\|\Psi_s\|_{\mathcal{B}(L_2,L_2(\mathcal{H}))} = \Big\| \summ_k
(id-\mathsf{E}_k) T_{4 \cdot 2^{-k}} \Delta_{k+s}
\Big\|_{\mathcal{B}(L_2, L_2(\mathcal{H}))} \le
\mathrm{c}_{n,\gamma} \hskip1pt 2^{-\gamma s/2}.$$
\end{shifttheo}

\subsubsection*{\textnormal{A.2.1.} The norm of $\Phi_s$}
\label{SSPhi}

%\renewcommand{\theequation}{A.2}
%\addtocounter{equation}{-1}

Define
\begin{eqnarray*}
\phi_{R_x}(w) & = & \frac{1}{|R_x|} \, 1_{R_x}(w), \\
\psi_{\widehat{Q}_y}(z) & = & \frac{1}{|\widehat{Q}_y|} \,
\sum_{j=2}^{2^n} 1_{Q_y}(z) - 1_{Q_j}(z).
\end{eqnarray*}
where $R_x$ is the only cube in $\Q_k$ containing $x$ and $Q_y$ is
the only cube in $\Q_{k+s}$ containing $y$. Moreover, the cubes
$Q_2, Q_3, \ldots, Q_{2^n}$ represent the remaining cubes in
$\Q_{k+s}$ sharing dyadic father with $Q_y$. Arguing as in
\cite{P2}, it is easily checked that the kernel $k_{s,k}(x,y)$ of
$\mathsf{E}_k T \Delta_{k+s}$ has the form
\begin{equation} \label{kernel}
k_{s,k}(x,y) = \left\langle T(\psi_{\widehat{Q}_y}), \phi_{R_x}
\right\rangle.
\end{equation}

\begin{Alemma} \label{prelimest}
If $T^*1=0$, the following estimates hold\hskip1pt$:$
\begin{itemize}
\item[a)] If $y \in \R^n \setminus 3 R_x$, we have
$$\big\| k_{s,k}(x,y) \big\|_\mathcal{H} \le \mathrm{c}_n \hskip1pt 2^{-
\gamma (k+s)} \frac{1}{|x-y|^{n+\gamma}}.$$
\item[b)] If $y \in 3 R_x \setminus R_x$, we have $$\big\|
k_{s,k}(x,y) \big\|_\mathcal{H} \le \mathrm{c}_{n,\gamma}
\hskip1pt 2^{- \gamma (k+s)} 2^{nk} \min \Bigg\{ \int_{R_x}
\frac{dw}{|w -\mathrm{c}_y|^{n+\gamma}}, s 2^{\gamma (k+s)}
\Bigg\}.$$
\item[c)] Similarly, if $y \in R_x$ we have $$\hskip14pt \big\|
k_{s,k}(x,y) \big\|_\mathcal{H} \le \mathrm{c}_{n,\gamma}
\hskip1pt 2^{- \gamma (k+s)} 2^{nk} \min \Bigg\{ \int_{\R^n
\setminus R_x} \frac{dw}{|w -\mathrm{c}_y|^{n+\gamma}}, s
2^{\gamma (k+s)} \Bigg\}.$$
\end{itemize}
The constant $\mathrm{c}_{n,\gamma}$ only depends on $n$ and
$\gamma;$ $\mathrm{c}_y$ denotes the center of the cube
$\widehat{Q}_y$.
\end{Alemma}

The main ingredients of the proof are the size/smoothness
estimates imposed on the kernel, plus the cancellation condition
\eqref{T*1=0} and the localization lemma given above. In
particular, the proof in \cite[Lemma 2.3]{P2} translates verbatim
to the Hilbert space valued context. Moreover, the result below is
a direct consequence of Lemma \ref{prelimest} and some
calculations provided in \cite{P2}.

\begin{Alemma} \label{prelimest2} Let us define
$$\begin{array}{rclcl}
\mathcal{S}^1_{s,k}(x) & = & \displaystyle \int_{\R^n} \big\|
k_{s,k}(x,y) \big\|_\mathcal{H} \, dy, \\ [10pt]
\mathcal{S}^2_{s,k}(y) & = & \displaystyle \int_{\R^n} \big\|
k_{s,k}(x,y) \big\|_\mathcal{H} \, dx.
\end{array}$$
Then, there exists a constant $\mathrm{c}_{n,\gamma}$ depending
only on $n,\gamma$ such that
$$\begin{array}{rclcl}
\mathcal{S}^1_{s,k}(x) & \le & \displaystyle
\frac{\mathrm{c}_{n,\gamma} \hskip1pt s}{2^{\gamma s}} &
\mbox{for all} \quad (x,k) \in \R^n \!\! \times \Z, \\
[10pt] \mathcal{S}^2_{s,k}(y) & \le & \hskip1pt
\mathrm{c}_{n,\gamma} \hskip1pt s & \mbox{for all} \quad \hskip1pt
(y,k) \in \R^n \!\! \times \Z.
\end{array}$$
\end{Alemma}

Now we are in position to complete our estimate for $\Phi_s$. In
fact, the argument we are giving greatly simplifies the one
provided in \cite{P2}. Namely, let us write $\Lambda_{s,k}$ for
$\mathsf{E}_k T \Delta_{k+s}$. Then, Lemma \ref{prelimest2} in
conjunction with Schur lemma provides the estimate
$\|\Lambda_{s,k}\|_{\mathcal{B}(L_2,L_2(\mathcal{H}))} \le
\mathrm{c}_{n,\gamma} s 2^{- \gamma s/2}$. On the other hand,
according to Cotlar lemma, it remains to check that $$\max \Big\{
\big\| \Lambda_{s,i}^* \Lambda_{s,j}^{\null}
\big\|_{\mathcal{B}(L_2)}, \big\| \Lambda_{s,i}^{\null}
\Lambda_{s,j}^* \big\|_{\mathcal{B}(L_2(\mathcal{H}))} \Big\} \,
\le \, \mathrm{c}_{n,\gamma} \hskip1pt s^2 \hskip1pt e^{- \gamma
s} \alpha_{i-j}^2$$ for some sumable sequence $(\alpha_k)_{k \in
\Z}$. Now, by the orthogonality of martingale differences, it
suffices to estimate the mappings $\Lambda_{s,i}^*
\Lambda_{s,j}^{\null}$ in $\mathcal{B}(L_2)$. Indeed, the only
nonzero mapping $\Lambda_{s,i}^{\null} \Lambda_{s,j}^*$ is the one
given by $i=j$ and $$\big\| \Lambda_{s,k}^{\null} \Lambda_{s,k}^*
\big\|_{\mathcal{B}(L_2(\mathcal{H}))} \le \big\| \Lambda_{s,k}
\big\|_{B(L_2, L_2(\mathcal{H}))}^2 \le \mathrm{c}_{n,\gamma} s^2
e^{-\gamma s},$$ by the estimate given above. To estimate the norm
of $\Lambda_{s,i}^* \Lambda_{s,j}^{\null}$, we assume (with no
loss of generality) that $i \ge j$. The martingale property then
gives $\mathsf{E}_i \mathsf{E}_j = \mathsf{E}_j$, so that
$\Lambda_{s,i}^* \Lambda_{s,j}^{\null} = \Delta_{i+s} T^*
\mathsf{E}_j T \Delta_{j+s}$. If we combine this with the estimate
deduced from Lemma \ref{prelimest2}, we get
\begin{eqnarray*}
\big\| \Lambda_{s,i}^* \Lambda_{s,j}^{\null}
\big\|_{\mathcal{B}(L_2)} & \le & \big\| \Lambda_{s+i-j,j}
\big\|_{B(L_2, L_2(\mathcal{H}))} \big\| \Lambda_{s,j}
\big\|_{B(L_2, L_2(\mathcal{H}))} \\ & \le & \mathrm{c}_{n,\gamma}
s (s + |i-j|) e^{-\gamma s} e^{- \gamma |i-j|/2}.
\end{eqnarray*}

\subsubsection*{\textnormal{A.2.2.} The norm of $\Psi_s$}
\label{SSPsi}

The arguments in this case follows a similar pattern to those used
for $\Phi_s$. The main differences are two. First, we can not use
\eqref{T*1=0} anymore but this is solved by the cancellation
produced by the term $(id - \mathsf{E}_k)$. Second, the
simplification with respect to the original argument in \cite{P2}
given above --using the martingale property-- does not work
anymore in this setting and we have to follow the complete
argument in \cite{P2} for this case. Nevertheless, the translation
of the original proof to the present setting is again
straightforward. We leave the details to the interested reader.

\subsection*{\textnormal{A.3.} The paraproduct argument}
\label{PA}

To complete the proof, we need to provide an alternative argument
for those Calder\'on-Zygmund operators failing the cancellation
condition \eqref{T*1=0}. Going back to our original decomposition
of $1_{\R^n \setminus \Sigma_{f,s}} T \! f$, we may write $1_{\R^n
\setminus \Sigma_{f,s}} T \! f = 1_{\R^n \setminus \Sigma_{f,s}}
(\Phi_s f + \Psi_s f)$. The second term $\Psi_s f$ is fine because
the quasi-orthogonal methods applied to it do not require
condition \eqref{T*1=0}. For the first term we need to use the
dyadic paraproduct associated to $\rho = T^*1$ $$\Pi_\rho(f) =
\sum_{j=-\infty}^\infty \Delta_j(\rho) \mathsf{E}_{j-1}(f).$$

\begin{Alemma} \label{LemmaHW}
If $T: L_2 \to L_2(\mathcal{H})$ is bounded, then so is
$\Pi_\rho$.
\end{Alemma}

\dem We have $$\|\Pi_\rho(f)\|_{L_2(\mathcal{H})}^2 =
\sum_{j=-\infty}^\infty \int_{\R^n} \|d_j \rho(x)\|_\mathcal{H}^2
|f_{j-1}(x)|^2 \, dx.$$ Since $\mathsf{E}_{j-1}(\|d_j
\rho(x)\|_\mathcal{H}) = \|d_j \rho(x)\|_\mathcal{H}$, we may
define $$R = \sum_{j=-\infty}^\infty \|d_j \rho(x)\|_\mathcal{H}
\, r_j$$ with $r_j$ the $j$-th Rademacher function. Then, it is
clear that $$\|\Pi_\rho(f)\|_{L_2(\mathcal{H})} = \|\Pi_R(f)\|_2
\le \|R\|_{\mathrm{BMO}_d(\R^n)} \|f\|_2.$$ On the other hand, it
is not difficult to check that the identities below hold
\begin{eqnarray*}
\|R\|_{\mathrm{BMO}_d(\R^n)} & = & \sup_{j \in \Z} \Big\|
\mathsf{E}_j \sum_{k > j} |dR_k|^2 \Big\|_\infty^\frac12 \\ & = &
\sup_{j \in \Z} \Big\| \mathsf{E}_j \sum_{k > j}
\|d\rho_k\|_\mathcal{H}^2 \Big\|_\infty^\frac12 \\ & = & \sup_{j
\in \Z} \, \sup_{Q \in \mathcal{Q}_j} \Big[ \frac{1}{|Q|} \int_Q
\big\| \rho(x) - \rho_Q \big\|_\mathcal{H}^2 \, dx \Big]^{\frac12} \\
& \sim & \sup_{j \in \Z} \, \sup_{Q \in \mathcal{Q}_j}
\inf_{\alpha_Q \in \mathcal{H}} \Big[ \frac{1}{|Q|} \int_Q \big\|
\rho(x) - \alpha_Q \big\|_\mathcal{H}^2 \, dx \Big]^{\frac12}.
\end{eqnarray*}
Then, we take $\alpha_Q = T^* 1_{\R^n \setminus 2Q}(c_Q)$ with
$c_Q$ the center of $Q$ and obtain
\begin{eqnarray*}
\lefteqn{\null \hskip-20pt \Big[ \frac{1}{|Q|} \int_Q \big\|
\rho(x) - \alpha_Q \big\|_\mathcal{H}^2 \, dx \Big]^{\frac12}}
\\ & \le & \Big[ \frac{1}{|Q|} \int_Q
\big\| T^*1_{2Q}(x) \big\|_\mathcal{H}^2 \, dx \Big]^{\frac12} \\
& + & \Big[ \frac{1}{|Q|} \int_Q \big\| T^* 1_{\R^n \setminus
2Q}(x) - T^* 1_{\R^n \setminus 2Q}(c_Q) \big\|_\mathcal{H}^2 dx
\Big]^{\frac12} = \mathsf{A} + \mathsf{B}.
\end{eqnarray*}
The fact that $T^*: L_2 \to L_2(\mathcal{H})$ is bounded yields
$$\mathsf{A} \le \Big[ \frac{1}{|Q|} \int_{\R^n} \big\|
T^*1_{2Q}(x) \big\|_\mathcal{H}^2 \, dx \Big]^\frac12 \lesssim
\frac{1}{\sqrt{|Q|}} \, \|1_{2Q}\|_2 \le 2^{\frac{n}{2}}.$$ On the
other hand, Lipschitz smoothness gives for $x \in Q$ $$\big\| T^*
1_{\R^n \setminus 2Q}(x) - T^* 1_{\R^n \setminus 2Q}(c_Q)
\big\|_\mathcal{H} \le \int_{\R^n \setminus 2Q} \big\| k^*(x,y) -
k^*(c_Q,y) \big\|_\mathcal{H} \, dy \le \mathrm{c}_n.$$ This
automatically gives an absolute bound for $\mathsf{B}$ and the
result follows. \fin

According to Lemma \ref{LemmaHW}, the dyadic paraproduct
$\Pi_\rho: L_2 \to L_2(\mathcal{H})$ defines a bounded operator.
Therefore, regarding its adjoint as a mapping $L_2 \to
L_2(\mathcal{H})$ via the identifications explained above, we find
a bounded map $$\Pi_\rho^*(f) = \sum_{j = - \infty}^\infty
\mathsf{E}_{j-1} \big( \langle \Delta_j(\rho), \cdot \rangle f
\big).$$ This allows us to consider the decomposition
$$T = T_0 + \Pi_{\rho}^*.$$ Now we go back to our estimate.
Following \cite{P2}, we have $$1_{\R^n \setminus \Sigma_{f,s}}
\summ_k \mathsf{E}_k \Pi_\rho^* \Delta_{k+s} f = 0.$$ Therefore,
it suffices to see that $T_0$ satisfies the following estimate
$$\Big\| \summ_k \mathsf{E}_k T_0 \Delta_{k+s}
\Big\|_{\mathcal{B}(L_2,L_2(\mathcal{H}))} \le
\mathrm{c}_{n,\gamma} s \, e^{- \gamma s/2}.$$ It is clear that
$T_0^*1=0$, so that condition \eqref{T*1=0} holds for $T_0$.
Moreover, according to our previous considerations, $T_0: L_2 \to
L_2(\mathcal{H})$ is bounded and its kernel satisfies the size
estimate since the same properties hold for $T$ and $\Pi_\rho^*$.
In particular, the only problem to apply the quasi-orthogonal
argument to $T_0$ is to verify that its kernel satisfies Lipschitz
smoothness estimates. We know by hypothesis that $T$ does.
However, the dyadic paraproduct only satisfies dyadic analogues.
It is nevertheless enough. Indeed, since Lemma \ref{prelimest2}
follows automatically from Lemma \ref{prelimest}, it suffices to
check the latter. However, following the argument in \cite{P2}, it
is easily seen that the instances in Lemma \ref{prelimest} where
the smoothness of the kernel is used equally work (giving rise to
$0$) in the dyadic setting.

\section*{Appendix B. Background on noncommutative
integration}

\renewcommand{\theequation}{B\arabic{equation}}

\setcounter{equation}{0}

We end this article with a brief survey on noncommutative $L_p$
spaces and related topics that have been used along the paper.
Most of these results are well-known to experts in the field. The
right framework for a noncommutative analog of the classical
measure theory and integrations is the theory of von Neumann
algebras. We refer to \cite{KR,Ta} for a systematic study of von
Neumann algebras and to the recent survey by Pisier/Xu \cite{PX2}
for a detailed exposition of noncommutative $L_p$ spaces.

\subsection*{\textnormal{B.1.} Noncommutative $L_p$}
\label{NCLp}

A \emph{von Neumann algebra} is a weak-operator closed
$\mathrm{C}^*$-algebra. By the Gelfand-Naimark-Segal theorem, any
von Neumann algebra $\M$ can be embedded in the algebra
$\mathcal{B}(\mathcal{H})$ of bounded linear operators on some
Hilbert space $\mathcal{H}$. In what follows we will identify $\M$
with a subalgebra of $\mathcal{B(H)}$. The positive cone $\M_+$ is
the set of positive operators in $\M$. A \emph{trace} $\tau: \M_+
\to [0,\infty]$ on $\M$ is a linear map satisfying the tracial
property $\tau(a^*a) = \tau(aa^*)$. A trace $\tau$ is
\emph{normal} if $\sup_\alpha \tau(a_\alpha) = \tau(\sup_\alpha
a_\alpha)$ for any bounded increasing net $(a_\alpha)$ in $\M_+$;
it is \emph{semifinite} if for any non-zero $a \in \M_+$, there
exists $0 < a' \le a$ such that $\tau(a') < \infty$ and it is
\emph{faithful} if $\tau(a) = 0$ implies $a = 0$. Taking into
account that $\tau$ plays the role of the integral in measure
theory, all these properties are quite familiar. A von Neumann
algebra $\M$ is called \emph{semifinite} whenever it admits a
normal semifinite faithful (\emph{n.s.f.} in short) trace $\tau$.
Recalling that any operator $a$ can be written as a linear
combination $a_1 - a_2 + ia_3 - ia_4$ of four positive operators,
we can extend $\tau$ to the whole algebra $\M$. Then, the tracial
property can be restated in the familiar way $\tau(ab) = \tau(ba)$
for all $a,b \in \M$.

According to the GNS construction, it is easily seen that the
noncommutative analogs of measurable sets (or equivalently
characteristic functions of those sets) are orthogonal
projections. Given $a \in \M_+$, the support projection of $a$ is
defined as the least projection $q$ in $\M$ such that $qa = a =
aq$ and will be denoted by $\mbox{supp} \hskip1pt a$. Let
$\mathcal{S}_+$ be the set of all $a \in \M_+$ such that
$\tau(\mbox{supp} \hskip1pt a) < \infty$ and set $\mathcal{S}$ to
be the linear span of $\mathcal{S}_+$. If we write $|x|$ for the
operator $(x^*x)^{\frac12}$, we can use the spectral measure
$\gamma_{|x|}: \R_+ \to \mathcal{B}(\mathcal{H})$ of the operator
$|x|$ to define
$$|x|^p = \int_{\R_+} s^p \, d \gamma_{|x|}(s) \quad
\mbox{for} \quad 0 < p < \infty.$$ We have $x \in \mathcal{S}
\Rightarrow |x|^p \in \mathcal{S}_+ \Rightarrow \tau(|x|^p) <
\infty$. If we set $\|x\|_p = \tau( |x|^p )^{\frac1p}$, it turns out
that $\| \ \|_p$ is a norm in $\mathcal{S}$ for $1 \le p < \infty$
and a $p$-norm for $0 < p < 1$. Using that $\mathcal{S}$ is a
$w^*$-dense $*$-subalgebra of $\M$, we define the
\emph{noncommutative $L_p$ space} $L_p(\M)$ associated to the pair
$(\M, \tau)$ as the completion of $(\mathcal{S}, \| \ \|_p)$. On the
other hand, we set $L_\infty(\M) = \M$ equipped with the operator
norm. Many fundamental properties of classical $L_p$ spaces, like
duality, real and complex interpolation... have been transferred to
this setting. The most important properties for our purposes are
\begin{itemize}
\item H\"older inequality. If $1/r = 1/p+1/q$, we have $\|ab\|_r
\le \|a\|_p \|b\|_q$.
\item The trace $\tau$ extends to a continuous functional on
$L_1(\M)$: $|\tau(x)| \le \|x\|_1$.
\end{itemize}
See \cite{PX2} for $L_p$ spaces over type III algebras. Let us
recall a few examples:

\begin{itemize}
\item[\textbf{(a)}] \textbf{Commutative $L_p$ spaces.} Let $\M$ be
commutative and semifinite. Then there exists a semifinite measure
space $(\Omega, \Sigma, \mu)$ for which $\M = L_\infty(\Omega)$
and $L_p(\M) = L_p(\Omega)$ with the \emph{n.s.f.} trace $\tau$
determined by
$$\tau(f) = \int_\Omega f(\omega) \, d \mu(\omega).$$

\item[\textbf{(b)}] \textbf{Semicommutative $L_p$ spaces.} Given a
measure space $(\Omega, \Sigma, \mu)$ and a semifinite von Neumann
algebra $(\mathcal{N}, \tau)$. We consider the von Neumann algebra
$$(\M,\mu\otimes\tau) = \Big( L_\infty(\Omega)\bar\otimes \mathcal{N},
\int_\Omega \tau ( \cdot) \, d \mu \Big).$$ In this case,
$$L_p(\M) = L_p \big( \Omega; L_p(\mathcal{N})
\big),$$ the $L_p$-space of $L_p(\mathcal{N})$-valued
Bochner-integrable functions, for all $p<\infty$.

\vskip5pt

\item[\textbf{(c)}] \textbf{Schatten $p$-classes.} Let $\M =
\mathcal{B}(\mathcal{H})$ with the standard trace
$$\mbox{tr}(x) = \summ_{\lambda} \langle x e_{\lambda}, e_{\lambda}
\rangle_\mathcal{H},$$ where $(e_{\lambda})_{\lambda}$ is any
orthonormal basis of $\mathcal{H}$. Then, the associated $L_p$
space is called the Schatten $p$-class
$\mathcal{S}_p(\mathcal{H})$. When $\mathcal{H}$ is separable, the
Schatten $p$-class is the noncommutative analog of $\ell_p$, which
embeds isometrically into the diagonal of $\mathcal{S}_p$.

\vskip5pt

\item[\textbf{(d)}] \textbf{Hyperfinite II$_1$ factor.} Let $M_2$
be the algebra of $2 \times 2$ matrices equipped with the
normalized trace $\sigma = \frac{1}{2} \mbox{tr}$. A description
of the so-called hyperfinite II$_1$ factor $\mathcal{R}$ is by the
following von Neumann algebra tensor product
$$(\mathcal{R},\tau) = \overline{\bigotimes_{n \ge 1}} (M_2,
\sigma).$$ That is, $\mathcal{R}$ is the von Neumann algebra
generated by all elementary tensors $x_1 \otimes \cdots \otimes
x_n \otimes \mathbf{1} \otimes \mathbf{1} \cdots$ and the trace
$\tau$ is the unique normalized trace on $\mathcal{R}$ which is
determined by $$\tau \Big( x_1 \otimes \cdots \otimes x_n \otimes
\mathbf{1} \otimes \mathbf{1} \cdots \Big) = \prod_{k=1}^n
\sigma(x_k).$$ $L_p(\mathcal{R})$ may be regarded as a
noncommutative analog of $L_p[0,1]$.
\end{itemize}

\noindent There are many other nice examples which we are
omitting, like free product von Neumann algebras, $q$-deformed
algebras, group von Neumann algebras... We refer to \cite{X} for a
more detailed explanation.

\begin{Bremark}
\emph{Assume that the Hilbert space ${\mathcal H}$ has a countable
orthonormal basis $(e_k)_{k \in \Z}$ and consider the associated
unit vectors $e_{m,n}$ of $\mathcal{B}({\mathcal H})$ which are
determined by the relation $e_{m,n}(x) = \langle x, e_n \rangle
e_m$. Let $T_{+}$ and $T_{-}$ be the projections from
$\mathcal{B}(\mathcal{H})$ onto the subspaces $$\Lambda_+ =
\mathrm{span} \big\{ e_{m,n} \, | \, m>n \big\} \quad \mbox{and}
\quad \Lambda_- = \mathrm{span} = \big\{ e_{m,n} \, | \, m<n
\big\}$$ respectively. $T_{+}$ and $T_{-}$ are the first examples
of noncommutative singular integrals and $T_{+}-T_{-}$ is a
nononcommutative analog of the classical Hilbert transform. In
fact, when ${\mathcal H} = L_2(\mathbb{T})$ is the space of all
$L_2$-integrable functions on the unit circle and $e_k =
\mathrm{exp}(i k \, \cdot)$, we embed $L_\infty (\Bbb{T})$ into
$\mathcal{B}({\mathcal H})$ via the map $$\Psi: f \in
L_\infty(\mathbb{T}) \mapsto  \Psi(f) \in \mathcal{B}({\mathcal
H}) \quad \mbox{with} \quad \Psi(f)[g] = \overline{f} g,$$ for any
$g \in L_2({\Bbb T})$. Then $\Psi(e_k) = \sum_m e_{m,m-k}$ and
\[ (T_{+}-T_{-})(\Psi (f))=-i\Psi (Hf), \]
for any $f \in L_2({\Bbb T}) \bigcap L_\infty({\Bbb T})$. Here $H$
denotes the classical Hilbert transform on ${\Bbb T}$.}
\end{Bremark}

\subsection*{\textnormal{B.2.} Noncommutative symmetric spaces}

Let $$\M' = \Big\{ b \in \mathcal{B}(\mathcal{H}) \ \big| \, ab =
ba \ \mbox{for all} \ a \in \M \Big\}$$ be the commutant of $\M$.
A closed densely-defined operator on $\mathcal{H}$ is
\emph{affiliated} with $\M$ when it commutes with every unitary
$u$ in the commutant $\M'$. Recall that $\M = \M''$ and this
implies that every $a \in \M$ is affiliated with $\M$. The
converse fails in general since we may find unbounded operators.
If $a$ is a densely defined self-adjoint operator on $\mathcal{H}$
and $a = \int_{\R} s \hskip1pt d \gamma_a(s)$ is its spectral
decomposition, the spectral projection $\int_{\mathcal{R}} d
\gamma_a(s)$ will be denoted by $\chi_{_\mathcal{R}}(a)$. An
operator $a$ affiliated with $\mathcal{M}$ is
\emph{$\tau$-measurable} if there exists $s > 0$ such that $$\tau
\big\{ |a| > s \big\}=\tau \big( \chi_{(s,\infty)} (|a|) \big) <
\infty.$$ The \emph{generalized singular-value} $\mu(a): \R_+ \to
\R_+$ is defined by
$$\mu_t (a) = \inf \Big\{ s
> 0 \, \big| \ \tau \big\{ |x| > s \big\} \le t \Big\}.$$
This provides us with a noncommutative analogue of the so-called
non-increasing rearrangement of a given function. We refer to
\cite{FK} for a detailed exposition of the function $\mu(a)$ and
the corresponding notion of convergence in measure.
If $L_0(\M)$ denotes the $*$-algebra of $\tau$-measurable
operators, we have the following equivalent definition of $L_p$
$$L_p(\M) = \Big\{a \in L_0(\M) \, \big| \ \Big( \int_{\R_+}
\mu_t(a)^p \, dt \Big)^{\frac1p} < \infty \Big\}.$$ The same
procedure applies to symmetric spaces. Given the pair $(\M,\tau)$,
let $\mathrm{X}$ be a rearrangement invariant quasi-Banach
function space on the interval $(0, \tau(\mathbf{1}_\M))$. The
\emph{noncommutative symmetric space} $\mathrm{X}(\mathcal{M})$ is
defined by
$$\mathrm{X}(\mathcal{M}) = \Big\{a \in L_0(\mathcal{M}) \,
\big| \ \mu(a) \in \mathrm{X} \Big\} \quad \text{with} \quad
\left\| a \right\|_{\mathrm{X}(\M)} = \|\mu(a)\|_{\mathrm{X}}.$$
It is known that $\mathrm{X}(\mathcal{M})$ is a Banach (resp.
quasi-Banach) space whenever $\mathrm{X}$ is a Banach (resp.
quasi-Banach) function space. We refer the reader to
\cite{DDdP,X1} for more in depth discussion of this construction.
Our interest in this paper is restricted to noncommutative
$L_p$-spaces and \emph{noncommutative weak $L_1$-spaces}.
Following the construction of symmetric spaces of measurable
operators, the noncommutative weak $L_1$-space
$L_{1,\infty}(\mathcal{M})$,  is defined as the set of all $a$ in
$L_0(\mathcal{M})$ for which the quasi-norm
\[ \left\|a\right\|_{1,\infty} = \sup_{t > 0} \, t \hskip1pt
\mu_t(x) = \sup_{\lambda > 0} \, \lambda \hskip1pt \tau \Big\{ |x|
> \lambda \Big\}\] is finite. As in the commutative case, the
noncommutative weak $L_1$ space satisfies a quasi-triangle
inequality that will be used below with no further reference.
Indeed, the following inequality holds for $a_1, a_2 \in
L_{1,\infty}(\M)$ $$\lambda \, \tau \Big\{ |a_1+a_2| > \lambda
\Big\} \le \lambda \, \tau \Big\{ |a_1| > \lambda/2 \Big\} +
\lambda \, \tau \Big\{ |a_2| > \lambda/2 \Big\}.$$

\subsection*{\textnormal{B.3.} Noncommutative martingales}

Consider a von Neumann subalgebra (a weak$^*$ closed
$*$-subalgebra) $\mathcal{N}$ of a semifinite von Neumann algebra
$(\mathcal{M},\tau)$. A \emph{conditional expectation}
$\mathcal{E}: \mathcal{M} \to \mathcal{N}$ is a positive
contractive projection from $\mathcal{M}$ onto $\mathcal{N}$. The
conditional expectation $\mathcal{E}$ is called \emph{normal} if
the adjoint map $\mathcal{E}^*$ satisfies
$\mathcal{E}^*(\mathcal{M}_*) \subset \mathcal{N}_*$. In this
case, there is a map $\mathcal{E}_*: \mathcal{M}_* \rightarrow
\mathcal{N}_*$ whose adjoint is $\mathcal{E}$. Such normal
conditional expectation exists if and only if the restriction of
$\tau$ to the von Neumann subalgebra $\mathcal{N}$ remains
semifinite, see e.g. Theorem 3.4 in \cite{Ta}. This is always the
case when $\tau(\mathbf{1}_\M)<\infty$. Any such conditional
expectation is trace preserving (i.e. $\tau \circ \mathcal{E} =
\tau$) and satisfies the bimodule property
\[ \mathcal{E}(a_1 b \hskip1pt a_2) = a_1 \mathcal{E}(b) \hskip1pt
a_2 \quad \mbox{for all} \quad a_1, a_2 \in \mathcal{N} \
\mbox{and} \ b \in \mathcal{M}.\]
Let $(\mathcal{M}_k)_{k \ge 1}$  be an increasing sequence of von
Neumann subalgebras of $\mathcal{M}$ such that the union of the
$\mathcal{M}_k$'s is weak$^*$ dense in $\mathcal{M}$. Assume that
for every $k \ge 1$, there is a normal conditional expectation
$\mathcal{E}_k: \mathcal{M} \to \mathcal{M}_k$. Note that for
every $1 \le p < \infty$ and $k \ge 1$, $\mathcal{E}_k$ extends to
a positive contraction $\mathcal{E}_k: L_p(\mathcal{M}) \to
L_p(\mathcal{M}_k)$. A \emph{noncommutative martingale} with
respect to the filtration $(\mathcal{M}_k)_{k \ge 1}$ is a
sequence $a = (a_k)_{k \ge 1}$ in $L_1(\mathcal{M})$ such that
\[
\mathcal{E}_j(a_k) = a_j \quad \mbox{for all} \quad 1 \le j \le k
< \infty.
\]
If additionally $a \subset L_p(\mathcal{M})$ for some $1 \le p \le
\infty$ and $\|a\|_p = \sup_{k \ge 1} \|a_k\|_p < \infty$, then
$a$ is called an \emph{$L_p$-bounded martingale}. Given a
martingale $a = (a_k)_{k \ge 1}$, we assume the convention that
$a_0 = 0$. Then, the martingale difference sequence $da =(da_k)_{k
\ge 1}$ associated to $x$ is defined by $da_k = a_k - a_{k-1}$.

\noindent Let us now comment some examples of noncommutative
martingales:
\begin{itemize}
\item[\textbf{(a)}] \textbf{Classical martingales.} Given a
commutative finite von Neumann algebra $\M$ equipped with a
normalized trace $\tau$ and a filtration $(\M_n)_{n \ge 1}$, there
exists a probability space $(\Omega, \Sigma, \mu)$ and an
increasing sequence $(\Sigma_n)_{n \ge 1}$ of $\sigma$-subalgebras
satisfying $$L_p(\M) = L_p(\Omega,\Sigma,\mu) \quad \mbox{and}
\quad L_p(\M_n) = L_p(\Omega, \Sigma_n, \mu).$$ Thus, classical
martingales are a form of noncommutative martingales. \vskip5pt
\item[\textbf{(b)}] \textbf{Semicommutative martingales.} Let
$(\Omega, \Sigma, \mu)$ be a probability space and $(\mathcal{N},
\tau)$ be a semifinite von Neumann algebra. Given $(\Sigma_n)_{n
\ge 1}$ an increasing filtration of $\sigma$-subalgebras of
$\Sigma$, we consider the filtration $$(\M_n,\tau) = \Big(
L_\infty(\Omega,\Sigma_n, \mu)\bar\otimes \mathcal{N}, \int_\Omega
\tau(\cdot) \, d \mu \Big),$$ of $$(\M,\tau) = \Big(
L_\infty(\Omega,\Sigma, \mu) \bar\otimes \mathcal{N}, \int_\Omega
\tau(\cdot) \, d \mu \Big).$$ In this case, conditional
expectations are given by $\mathsf{E}_n = \mathbb{E}_n \otimes
id_{\mathcal{N}}$ where $\mathbb{E}_n$ denotes the conditional
expectation $(\Omega, \Sigma) \to (\Omega,\Sigma_n)$. Once more
(in this particular setting) noncommutative martingales can be
viewed as vector valued commutative martingales. \vskip5pt
\item[\textbf{(c)}] \textbf{Finite martingales.} When dealing with
Schatten $p$-classes, there is no natural finite trace unless we
work in the finite-dimensional case $\mathcal{S}_p(n)$ where we
consider the normalized trace $\sigma = \frac1n \mbox{tr}$. A
natural filtration is obtained taking $(\M_k,\sigma)$ to be the
subalgebra of $k \times k$ matrices (i.e. with vanishing entries
in the last $n-k$ rows and columns). This choice is useful to
obtain certain counterexamples, see \cite{JX3}. \vskip5pt
\item[\textbf{(d)}] \textbf{Dyadic martingales.} Let
$\varepsilon_1, \varepsilon_2, \varepsilon_3 \ldots$ be a
collection of independent $\pm 1$ Bernoullis. Classical dyadic
martingales are constructed over the filtration $\Sigma_n =
\sigma(\varepsilon_1, \varepsilon_2, \ldots, \varepsilon_n)$. The
noncommutative analog consists of a filtration $(\mathcal{R}_n)_n$
in the hyperfinite $\mathrm{II}_1$ factor as follows
$$(\mathcal{R}_n,\tau) = \overline{\bigotimes_{1 \le m \le n}} \,
(M_2,\sigma).$$ $\mathcal{R}_n$ embeds into $\mathcal{R}$ by means
of $a_1 \otimes \cdots \otimes a_n \mapsto a_1 \otimes \cdots
\otimes a_n \otimes \mathbf{1} \otimes \mathbf{1} \cdots$.
Moreover, given $1 \le n \le r$, the conditional expectation
$\mathsf{E}_n: \mathcal{R} \to \mathcal{R}_n$ is determined by
$$\null \hskip30pt \mathsf{E}_n \big( a_1 \otimes \cdots \otimes
a_r \otimes \mathbf{1} \otimes \mathbf{1} \cdots \big) = \Big(
\prod_{k=n+1}^r \sigma(a_k) \Big) a_1 \otimes \cdots \otimes a_n
\otimes \mathbf{1} \otimes \mathbf{1} \cdots$$
\end{itemize}

As we did with noncommutative $L_p$ spaces, we omit some standard
examples like free martingales, $q$-deformed martingales or
noncommutative martingales on the group algebra of a discrete
group. We refer again to \cite{X} for a more in depth exposition.

The theory of noncommutative martingales has achieved considerable
progress in recent years. The renewed interest on this topic
started from the fundamental paper of Pisier and Xu \cite{PX1},
where they introduced a new functional analytic approach to study
Hardy spaces and the Burkholder-Gundy inequalities for
noncommutative martingales. Shortly after, many classical
inequalities have been transferred to the noncommutative setting.
A noncommutative analogue of Doob's maximal function \cite{J1},
the noncommutative John-Nirenberg theorem \cite{JM}, extensions of
Burkholder inequalities for conditioned square functions
\cite{JX1} and related weak type inequalities \cite{R1,R3}; see
\cite{PR} for a simpler approach to some of them.

\vskip3pt

\noindent \textbf{Acknowledgements.} Tao Mei was partially
supported by a Young Investigator Award of the National Science
Foundation supported summer workshop in Texas A\&M University
2007. Javier Parcet has been partially supported by \lq Programa
Ram{\'o}n y Cajal, 2005\rq${}$ and by Grants MTM2007-60952,
CCG07-UAM/ESP-1664 and CCG08-CSIC/ESP-3485, Spain.

\bibliographystyle{amsplain}

\

\

\hfill \noindent \textbf{Tao Mei} \\
\null \hfill Department of Mathematics
\\ \null \hfill University of Illinois at Urbana-Champaign \\
\null \hfill 1409 W. Green St. Urbana, IL 61891. USA \\
\null \hfill\texttt{mei@math.uiuc.edu}

\

\hfill \noindent \textbf{Javier Parcet} \\
\null \hfill Instituto de Ciencias Matem{\'a}ticas \\ \null \hfill
CSIC-UAM-UC3M-UCM \\ \null \hfill Consejo Superior de
Investigaciones Cient{\'\i}ficas \\ \null \hfill Serrano 121.
28006, Madrid. Spain
\\ \null \hfill\texttt{javier.parcet@uam.es}
\end{document}